\documentclass[11pt]{amsart}
\usepackage{amsmath,amssymb,amsthm,mathtools}
\usepackage[hidelinks]{hyperref}
\usepackage[margin=0.82in]{geometry}
\newtheorem{theorem}{Theorem}[section]
\newtheorem{proposition}[theorem]{Proposition}
\newtheorem{lemma}[theorem]{Lemma}
\newtheorem{corollary}[theorem]{Corollary}
\newtheorem{example}[theorem]{Example}
\newtheorem*{fixedexpectationlemma}{Fixed-expectation lemma}
\theoremstyle{remark}
\newtheorem{remark}[theorem]{Remark}

\newcommand{\C}{\mathbb C}
\newcommand{\M}{\mathcal M_q}
\newcommand{\spn}{\operatorname{span}}
\newcommand{\Int}{\operatorname{int}}
\newcommand{\tr}{\operatorname{tr}}
\newcommand{\F}{\mathcal F}
\newcommand{\PSD}{\mathcal D_n}
\title[Multiplicities in the $q$-numerical range]{Multiplicities of Points in the $q$-Numerical Range}
\author{Tin-Yau Tam}
\address{Department of Mathematics and Statistics, University of Nevada, Reno, Nevada, USA}
\email{ttam@unr.edu}

\author{Yayi Yuan}
\address{College of Mathematics Science, Inner Mongolia Normal University, Hohhot, Inner Mongolia,China}
\email{yuanyayi@mails.imnu.edu.cn}

\thanks{Corresponding author: Yayi Yuan
(\texttt{yuanyayi@mails.imnu.edu.cn}).}

\date{November 24, 2026}
\begin{document}
\begin{abstract}
For $A\in M_n(\C)$, $n\ge2$, and $|q|<1$, we study two
multiplicities of points of the $q$-numerical range $W_q(A)$: the complex dimension of the span of its
generating first vectors and the real dimension of the span of their
rank-one orthogonal projections. For $0<|q|<1$, the first multiplicity
is $n$ at every relative interior point, and the second is $n^2$ when
$n\ge3$. At boundary points of nonscalar matrices, an optimal scalar shift
and compression to the top singular subspace yield exact multiplicity
formulas. In
particular, the projection multiplicity is $m^2-d$, where $m$ is the
generating-vector multiplicity and $d$ is the affine dimension of the
minimal face containing the prescribed expectation in the compressed
numerical range. We obtain explicit results for Hermitian and normal
matrices and positive weighted shifts, including a classification in
dimension two and parameter-dependent transitions. The endpoint $q=0$
is treated separately, where a two-dimensional interior exception occurs.
\end{abstract}
\maketitle

\section{Introduction}
Let $A\in M_n(\C)$ with $n\ge2$, and let $q\in\C$ satisfy
$|q|\le1$. The $q$-numerical range is
\[
W_q(A)=\{x^*Ay:\ x,y\in\C^n,\ \|x\|=\|y\|=1,\ x^*y=q\}.
\]
For $\xi\in W_q(A)$, define the generating-vector set and its
multiplicity by
\begin{align*}
\M(A,\xi)&=\{x\in\C^n:\ \|x\|=1 \ \text{and there exists}\ y \ \text{such that} \ \|y\|=1,\ x^*y=q,\ x^*Ay=\xi\},\\
m_q(A,\xi)&=\dim_{\C}\spn\M(A,\xi).
\end{align*}
Alongside this complex-linear multiplicity, we also consider the real-linear
span of the associated rank-one orthogonal projections,
\begin{equation}\label{eq:projection-multiplicity}
 p_q(A,\xi)=\dim_{\mathbb R}\operatorname{span}_{\mathbb R}
       \{xx^*:x\in\mathcal M_q(A,\xi)\}.
\end{equation}
The first multiplicity measures how many linearly independent vectors
can occur as the first member of a pair generating $\xi$. The second
records the linear information carried by the corresponding orthogonal
projections, which do not change when a generating vector is multiplied
by a scalar of modulus one. If $V=\spn\mathcal M_q(A,\xi)$ has
complex dimension $m$, then every such projection is supported on $V$.
Since the real vector space of Hermitian operators on $V$ has dimension
$m^2$, we have
\[
  1\le p_q(A,\xi)\le m_q(A,\xi)^2\le n^2.
\]
The two multiplicities therefore measure different features of the
same generating fiber. In particular, the equality $m_q(A,\xi)=n$
does not imply $p_q(A,\xi)=n^2$; the boundary examples below make
this distinction explicit.

For a convex set $K\subseteq\C$, denote its affine hull by
$\operatorname{aff}K$, its relative interior in $\operatorname{aff}K$
by $\operatorname{ri}K$, and its ordinary planar interior by
$\operatorname{int}K$. In particular, the relative interior of a
segment excludes its endpoints, while that of a singleton is the
singleton itself.

For $q=\rho e^{i\phi}$, $0<\rho\le1$, replacing $y$ by
$e^{-i\phi}y$ gives
\[
W_q(A)=e^{i\phi}W_\rho(A),\qquad
\mathcal M_q(A,\xi)=\mathcal M_\rho(A,e^{-i\phi}\xi).
\]
Consequently,
\[
m_q(A,\xi)=m_\rho(A,e^{-i\phi}\xi),\qquad
p_q(A,\xi)=p_\rho(A,e^{-i\phi}\xi).
\]
The phase identity for the range is standard; see
Li, Mehta, and Rodman~\cite[Proposition~2.2(c)]{LiMehtaRodman}
and Li and Nakazato~\cite[(1.1)(d)]{LiNakazato1998}.
Thus it suffices to study real $0<q<1$ and the endpoint $q=0$.
Convexity of $W_q(A)$ was proved by Tsing~\cite{Tsing}; at $q=1$
it reduces to the Toeplitz--Hausdorff theorem. We count the first
generating vectors; the second vectors have the same multiplicity
(Remark~\ref{rem:equal-vector-multiplicity}). 

The corresponding problem for the classical numerical range was studied
by Mary R.~Embry~\cite[Theorem~1]{Embry}. When $q=1$, the constraint
$x^*y=1$ forces $y=x$, so $W_1(A)=W(A)$ and our definitions give
\[
\mathcal M_1(A,\xi)=\{x\in\C^n:\|x\|=1,\ x^*Ax=\xi\},\qquad
m_1(A,\xi)=\dim_{\C}\spn\mathcal M_1(A,\xi).
\]
Write $V_1(A,\xi)=\spn\mathcal M_1(A,\xi)$ and
$\mathcal C_A(\xi)=\{v\in\C^n:v^*Av=\xi\|v\|^2\}$.
The support function is
\[
 h_{W(A)}(\theta)=\max_{\eta\in W(A)}\operatorname{Re}(e^{-i\theta}\eta).
\]
In this notation, Embry's theorem asserts the following:
\begin{enumerate}
\item $\xi$ is an extreme point of $W(A)$ if and only if
$\mathcal C_A(\xi)$ is a complex linear subspace (and hence equals
$V_1(A,\xi)$).
\item If $\xi$ is a nonextreme boundary point and $L$ is a supporting
line at $\xi$ with outward normal $e^{i\theta}$, then
\[
 V_1(A,\xi)=\ker\bigl(\operatorname{Re}(e^{-i\theta}A)
      -h_{W(A)}(\theta)I\bigr).
\]
Thus $m_1(A,\xi)=\dim_{\C}\ker(\operatorname{Re}(e^{-i\theta}A)-h_{W(A)}(\theta)I)$.
\item If $W(A)$ has nonempty planar interior, then
$m_1(A,\xi)=n$ if and only if $\xi\in\operatorname{int}W(A)$.
\end{enumerate}
Embry's result provides the starting point for our investigation,
but the passage from $q=1$ to $0<q<1$ changes the problem. At $q=1$,
every generating pair has $x=y$, and a supporting line leads to an
eigenspace of a Hermitian part of $A$. For $0<q<1$, the two vectors
are distinct. A supporting pair instead arises from a top singular
subspace of a scalar shift of $A$, together with an additional
numerical-range constraint on that subspace. Thus Embry's
supporting-eigenspace description does not directly give the
multiplicity at a boundary point of $W_q(A)$.

\smallskip
\noindent\emph{Earlier results on $W_q(A)$.}
Tsing~\cite{Tsing} established convexity of $W_q(A)$.
Li, Mehta, and Rodman~\cite{LiMehtaRodman} developed basic properties,
including boundary regularity, and Li and Nakazato~\cite{LiNakazato1998}
studied its geometry and $q$-numerical radii. The description through
the Davis--Wielandt shell and the strict convexity of the range are
associated with Chien and Nakazato~\cite{ChienNakazato2002};
see also~\cite{LiMehtaRodman}. For normal matrices we refer to
Li~\cite{Li1998} and Chien and Nakazato~\cite{ChienNakazato2005},
and for weighted shifts to~\cite{LiNakazato1998,ChienNakazato,
ChienNakazatoPerturbation}. These works concern, in particular,
the location, shape, boundary, and radius of $W_q(A)$. We use these
known properties but ask a different question: for a 
\emph{fixed point} $\xi\in W_q(A)$, what are the dimensions of the
linear spaces generated by its representing vectors and projections?
The related multiplicity problem for joint numerical ranges was
studied by Cheung, Liu, and Tam~\cite{CheungLiuTam}.

\smallskip
\noindent\emph{Main results.}
For $0<q<1$, every relative interior point has full
vector multiplicity (Theorem~\ref{thm:all}):
\[
 m_q(A,\xi)=n\qquad
 (n\ge2,\ \xi\in\operatorname{ri}W_q(A)).
\]
For $n\ge3$, its generating projections also span the full real
space of Hermitian matrices, giving
$p_q(A,\xi)=n^2$ (Theorem~\ref{thm:interior-projection}).
The proof of the vector result uses a rank-one semidefinite
lifting; dimension two requires a separate argument because a fixed
first vector generates a circle rather than a disk.

For boundary points of a nonscalar matrix, we obtain an optimal
scalar-shift support formula and an exact reduction to the numerical
range of an operator compressed to a top singular subspace
(Theorems~\ref{thm:support-dual} and~\ref{thm:boundary-exact}).
A second Hermitian compression turns this into an eigenspace
formula (Theorem~\ref{thm:boundary-hermitian}). In particular,
boundary multiplicity need not be one, and full vector multiplicity
may occur on the boundary. For projections we obtain the
minimal-face formula (Theorem~\ref{thm:boundary-face-projection})
\[
 p_q(A,\xi)=m_q(A,\xi)^2-d,
 \qquad
 d=\dim_{\mathbb R}\operatorname{aff}F_q,
\]
where $F_q$ is the smallest face containing $q$ of the compressed
numerical range arising from the boundary reduction. Consequently,
projection multiplicity records information not contained in
vector multiplicity alone.

For example, let $A=\operatorname{diag}(1,\omega,\omega^2)$,
where $\omega=e^{2\pi i/3}$, and take $0<q<1/2$. Then $W_q(A)$
is the closed unit disk. Every interior point has multiplicities
$m_q=3$ and $p_q=9$, whereas every boundary point has
$m_q=3$ and $p_q=7$ (Examples~\ref{ex:roots-unity}
and~\ref{ex:root-projections}). Thus the vector multiplicity is
unchanged, even though the projection multiplicity detects the
boundary. We also obtain explicit classifications for Hermitian and
normal matrices (Section~\ref{sec:hermitian}) and positive weighted
nilpotent shifts (Section~\ref{sec:weighted}), including changes of
multiplicity as $q$ varies. The projection multiplicities are treated in
Section~\ref{sec:proj}; Section~\ref{sec:qzero} treats the endpoint $q=0$,
and Section~\ref{sec:complexq} reduces nonzero complex parameters
to positive real parameters.

\begin{remark}[Either generating vector gives the same multiplicity]
\label{rem:equal-vector-multiplicity}
For a generating pair $(x,y)$ of $\xi$, we have $\|x\|=\|y\|=1$,
$x^*y=q$, and $x^*Ay=\xi$. For every $A\in M_n(\C)$,
$|q|\le1$, and $\xi\in W_q(A)$,
\[
\dim_{\C}\spn\{x:\ (x,y)\ \text{generates }\xi\ \text{for some }y\}
=\dim_{\C}\spn\{y:\ (x,y)\ \text{generates }\xi\ \text{for some }x\}.
\]
Swapping the two vectors replaces $(A,q,\xi)$ by
$(A^*,\overline q,\overline\xi)$. For real $0<q<1$, equality of the
span dimensions follows from Theorems~\ref{thm:all} and
\ref{thm:boundary-exact}, and for $q=0$ from
Theorems~\ref{thm:q0interior}--\ref{thm:q0boundary}.
For nonreal $q$ use the phase identity above. When $|q|=1$,
$y=qx$. Equal dimensions do not imply equal subspaces.
\end{remark}

Throughout Sections~\ref{sec:prelim}--\ref{sec:proj} we assume
$0<q<1$ and write $s=\sqrt{1-q^2}$ unless stated otherwise.

\section{Geometry of generating pairs}\label{sec:prelim}
For a unit vector $x$, the admissible second vectors are
\[
y=qx+su,\qquad u\perp x,\quad\|u\|=1.
\]
For a fixed unit $x$, let
\[
\mu(x)=x^*Ax,\qquad
r_A(x)=\sqrt{\|A^*x\|^2-|x^*Ax|^2}.
\]
With $x$ fixed, the possible values of $x^*Ay$ form a disk of center $q\mu(x)$
and radius $sr_A(x)$ when $n\ge3$, and generally its boundary circle
when $n=2$. To see this, note that
\[
x^*Au=(A^*x)^*u=\bigl(P_{x^\perp}A^*x\bigr)^*u,
\qquad \|P_{x^\perp}A^*x\|=r_A(x).
\]
If $n=2$, the space $x^\perp$ is one-dimensional: as $u$ runs
through its unit vectors, the value $x^*Au$ runs through a circle.
If $n\ge3$, we can vary the component of $u$ along
$P_{x^\perp}A^*x$ and use an orthogonal component to keep
$\|u\|=1$. Hence every point in the disk occurs. The distinction
is important when we study the generating vectors, rather than
only the set $W_q(A)$.
The resulting union-of-disks description is known from Tsing's
representation and its Davis--Wielandt-shell formulation;
see~\cite{Tsing,LiNakazato1998,ChienNakazato2002}.
We recall it with the fixed-first-vector distinction needed below.
Define, for all $n\ge2$, the \emph{auxiliary disk}
\begin{equation}\label{eq:disk}
	D_q(A,x)=\{\xi\in\C:|\xi-q\mu(x)|\le sr_A(x)\}.
\end{equation}
\begin{lemma}[Disk-union representation]\label{lem:disk-union}
For every $n\ge2$ one has
\[ W_q(A)=\bigcup_{\|x\|=1}D_q(A,x). \]
For $n=2$ this is an equality only \emph{after taking the union}:
an individual disk need not consist of values generated with its
specified first vector.
\end{lemma}
\begin{proof}
For $n\ge3$ the equality follows directly from the preceding
sphere-to-disk description. For $n=2$, the fixed-first-vector circle
is contained in the convex set $W_q(A)$ by Tsing's theorem; hence its
convex hull, which is $D_q(A,x)$, is also contained in $W_q(A)$.
Conversely, any generated value lies on a corresponding circle.
Thus in dimension two the disk-union identity describes the range
correctly, but the membership $\xi\in D_q(A,x)$ alone does not imply
that this particular $x$ generates $\xi$.
\end{proof}
Li, Mehta, and Rodman~\cite[Theorem~2.7]{LiMehtaRodman}
proved the stronger inclusion
$q\operatorname{conv}\sigma(A)\subset\operatorname{int}W_q(A)$
for nonscalar $A$ and $0<q<1$. We retain the following short proof
of nonempty interior for use with the disk representation.
\begin{lemma}[Compactness and nonempty interior]\label{lem:compact-int}
Let $A\in M_n(\C)$, $n\ge2$, and $0<q<1$. Then $W_q(A)$ is
nonempty and compact. If $A$ is nonscalar, $W_q(A)$ has nonempty
interior in $\C$.
\end{lemma}
\begin{proof}
Consider the set of admissible generating pairs
\[
\mathcal P_q=\{(x,y)\in\C^n\times\C^n:\
 \|x\|=\|y\|=1,\ x^*y=q\}.
\]
It is a closed subset of the product of two unit spheres, and is
therefore compact. It is nonempty: for the standard basis vectors
$e_1,e_2$, the pair $x=e_1$, $y=qe_1+se_2$ is admissible.
The map $(x,y)\mapsto x^*Ay$ is continuous, so its image
$W_q(A)$ is compact and nonempty.

Suppose now that $A$ is nonscalar. There is a unit vector $x_0$
for which $A^*x_0$ is not a multiple of $x_0$. Otherwise every
nonzero vector would be an eigenvector of $A^*$; applying this to
the standard basis vectors and their pairwise sums shows that
$A^*$, hence $A$, is scalar. Thus $r_A(x_0)>0$.
By Lemma~\ref{lem:disk-union}, the disk $D_q(A,x_0)$ is contained
in $W_q(A)$, including when $n=2$. Its radius is
$sr_A(x_0)>0$, since $0<q<1$. Hence $W_q(A)$ contains an open
disk and has nonempty interior.
\end{proof}

\begin{lemma}[Strict-disk criterion]\label{lem:strict}
If $n\ge3$ and some unit vector $x_0\in \C^n$ satisfies
\begin{equation}\label{eq:strict}
 |\xi-qx_0^*Ax_0|<s\sqrt{\|A^*x_0\|^2-|x_0^*Ax_0|^2},
\end{equation}
then $m_q(A,\xi)=n$.
\end{lemma}
\begin{proof}
The strict inequality persists on an open neighborhood of $x_0$
in the unit sphere. When $n\ge3$ the full fixed-first-vector disks
are genuine, so all the nearby vectors generate $\xi$. A nonempty
open subset of the unit sphere spans $\C^n$.
\end{proof}
\begin{proposition}[Covariance]\label{prop:inv}
For unitary $U$, 
\[
m_q(U^*AU,\xi)=m_q(A,\xi).
\]
 Also, for
$\alpha\in\C\setminus\{0\}$ and $\beta\in\C$,
\[ m_q(\alpha A+\beta I,\alpha \xi+q\beta)=m_q(A,\xi).\]
\end{proposition}
\begin{proof}
Transform the admissible pair $(x,y)$ by $U$, or use
$x^*(\alpha A+\beta I)y=\alpha x^*Ay+q\beta$.
\end{proof}

\section{A rank-one semidefinite lifting and interior multiplicity}\label{sec:general}
To establish interior multiplicity, we shall show that each interior
value lies strictly inside one of the disks in
Lemma~\ref{lem:disk-union}. We use density matrices and a rank-one
reduction. Put
\[
 \PSD=\{X\in M_n(\C):X=X^*\ge0,\ \tr X=1\}.
\]
By spectral decomposition and convexity of the classical numerical
range, $\tr(BX)\in W(B)$ for every $X\in\PSD$. We also use the
finite-dimensional fact that compact convex sets are the convex
hulls of their extreme points.

\begin{lemma}\label{lem:purify}
Let $A\in M_n(\C)$, and suppose $\mu\in W(A)$, the classical
numerical range. Then
\begin{equation}\label{eq:rank-one}
\max_{\substack{X\in\PSD\\\tr(AX)=\mu}}
       \tr(AA^*X)
=\max_{\substack{\|x\|=1\\x^*Ax=\mu}}\|A^*x\|^2.
\end{equation}
In particular the maximum on the left has a rank-one optimizer.
\end{lemma}
\begin{proof}
The feasible set is nonempty, compact, and convex, subject to 
\[
\tr X=1, \quad \operatorname{Re}\tr(AX)=\operatorname{Re}\mu,\quad 
\operatorname{Im}\tr(AX)=\operatorname{Im}\mu.
\]
The maximum is attained at an extreme point. If such a point $X$
has rank $r\ge2$, the real space of Hermitian matrices supported on
$\operatorname{ran}X$ has dimension $r^2\ge4$. Hence there is a
nonzero supported Hermitian $K$ such that
\[
\tr K=0,\qquad\operatorname{Re}\tr(AK)=0,
\qquad\operatorname{Im}\tr(AK)=0.
\]
The three real equations leave a nonzero solution $K$.
Since $X$ is positive definite on its range, $X\pm\varepsilon K$
remain positive semidefinite for sufficiently small $\varepsilon>0$.
They are distinct feasible matrices with midpoint $X$, a
contradiction. Thus every extreme point has rank one.
\end{proof}

For $n\ge2$, introduce the function, with $X\in\PSD$ and $\xi\in\C$,
\begin{equation}\label{eq:liftF}
\F_A(X,\xi)=s^2\tr(AA^*X)
   -|\tr(AX)-q\xi|^2-s^2|\xi|^2.
\end{equation}
The function $\F_A$ is jointly concave in $(X,\xi)$. For $X=xx^*$,
\begin{equation}\label{eq:lift-identity}
 \F_A(xx^*,\xi)=s^2r_A(x)^2-|\xi-qx^*Ax|^2.
\end{equation}
Thus $\F_A(xx^*,\xi)\ge0$ if and only if
$\xi\in D_q(A,x)$.
\begin{proposition}\label{prop:lift}
For every $A\in M_n(\C)$, $n\ge2$,
\begin{equation}\label{eq:SDP-characterization}
 W_q(A)=\{\xi\in\C:\ \text{there exists}\ X\in\PSD  \ \text{such that}\
       \F_A(X,\xi)\ge0\}.
\end{equation}
If the inequality is strict for some $X$, it is strict for a rank-one $X$.
\end{proposition}
\begin{proof}
If $\xi\in W_q(A)$, its generating pair has a vector $x$ satisfying
$\F_A(xx^*,\xi)\ge0$ by~\eqref{eq:lift-identity}.
Conversely, given $X\in\PSD$, put $\mu=\tr(AX)$.
Since $\mu\in W(A)$, Lemma~\ref{lem:purify} gives a unit $x$
with $x^*Ax=\mu$ and $\|A^*x\|^2\ge\tr(AA^*X)$. Hence
\[\F_A(xx^*,\xi)\ge\F_A(X,\xi).\]
The disk-union Lemma~\ref{lem:disk-union} proves~\eqref{eq:SDP-characterization}.
The same inequality proves the assertion about strict feasibility.
\end{proof}
\begin{remark}
The set $\{(X,\xi):X\in\PSD,\ \F_A(X,\xi)\ge0\}$ is convex;
its projection onto the second coordinate is $W_q(A)$.
\end{remark}

\begin{theorem}\label{thm:strict-all}
Let $n\ge2$, $0<q<1$, and let $A$ be nonscalar. For every
$\xi\in\Int W_q(A)$ there exists a unit vector $x\in \C^n$ such that
\begin{equation}\label{eq:strict-univ}
             |\xi-qx^*Ax|<sr_A(x).
\end{equation}
This assertion includes $n=2$, although for $n=2$ it does not
say that $x$ itself generates $\xi$.
\end{theorem}
\begin{proof}
Since $A^*$ is nonscalar, there is a unit $x_0$ for which
$A^*x_0$ is not proportional to $x_0$; thus $r_A(x_0)>0$.
Set $\xi_0=qx_0^*Ax_0$. Then
$\F_A(x_0x_0^*,\xi_0)=s^2r_A(x_0)^2>0$.
The disk-union lemma shows that the positive-radius disk
$D_q(A,x_0)$ belongs to $W_q(A)$, even when $n=2$.
Hence $\xi_0$ is an interior point.

Given $\xi\in\Int W_q(A)$ with $\xi\ne \xi_0$, extend the ray from
$\xi_0$ through $\xi$ to a boundary point $\xi_1\in W_q(A)$, obtaining
$\xi=(1-t)\xi_0+t\xi_1$ for $0<t<1$.
Choose the first vector $x_1$ of a generating pair for $\xi_1$
and put
$X_t=(1-t)x_0x_0^*+tx_1x_1^*$.
By joint concavity,
\[
\F_A(X_t,\xi)\ge(1-t)\F_A(x_0x_0^*,\xi_0)
                 +t\F_A(x_1x_1^*,\xi_1)>0.
\]
Proposition~\ref{prop:lift} produces a unit $x$ with
$\F_A(xx^*,\xi)>0$, equivalent to~\eqref{eq:strict-univ}.
The case $\xi=\xi_0$ is immediate.
\end{proof}

\begin{theorem}\label{thm:all}
If $n\ge2$, $A\in M_n(\C)$, and $0<q<1$, then
\[
           m_q(A,\xi)=n\quad\text{for all }\xi\in\operatorname{ri}W_q(A).
\]
\end{theorem}
\begin{proof}
If $A=\lambda I$, then $W_q(A)=\{q\lambda\}$ and every unit
vector occurs in a generating pair, so the result is immediate.
For nonscalar $A$, Lemma~\ref{lem:compact-int} gives
$\operatorname{ri}W_q(A)=\operatorname{int}W_q(A)$.

For $n\ge3$, Theorem~\ref{thm:strict-all} and
Lemma~\ref{lem:strict} give $m_q(A,\xi)=n$.
For $n=2$, consider the continuous function on
$\C P^1\cong S^2$,
\[
 g_\xi([x])=s^2r_A(x)^2-|\xi-qx^*Ax|^2,
 \qquad\|x\|=1.
\]
The circle condition in dimension two is precisely
\begin{equation}\label{eq:zeros}
 [x]\text{ generates }\xi\quad\Longleftrightarrow\quad g_\xi([x])=0.
\end{equation}
By Theorem~\ref{thm:strict-all}, $g_\xi$ takes a positive value.
To find a negative value, distinguish the two possible spectra of $A^*$.

If $A^*$ has two distinct eigenvalues $\bar\lambda_1$ and
$\bar\lambda_2$, corresponding unit eigenvectors satisfy
$r_A(x_j)=0$, hence
$g_\xi([x_j])=-|\xi-q\lambda_j|^2$. Since $q>0$ and
$\lambda_1\ne\lambda_2$, at least one value is negative.
Otherwise, since $A$ is nonscalar, its repeated eigenvalue
$\lambda$ is defective, and a unitary Schur form gives
$A=\lambda I+cJ_2$ with $c\ne0$ (after unitary similarity).
The vector $e_2$ satisfies $r_A(e_2)=0$, so $g_\xi([e_2])<0$
unless $\xi=q\lambda$. At that remaining value, for a unit
$x=(a,b)^\top$ one computes
\[
 g_{q\lambda}([x])=|c|^2|a|^2(|a|^2-q^2).
\]
Taking $0<|a|^2<q^2$ makes this negative.

If $m_q(A,\xi)=1$, the zero set of $g_\xi$ on $\C P^1$ is a
singleton. Its complement is path connected, contradicting the
existence of both positive and negative values of $g_\xi$.
Hence $m_q(A,\xi)=2$. 
\end{proof}

\begin{remark}[Why $q>0$ matters in dimension two]\label{rem:q0}
For $A=J_2$ and $q=0$ one has $W_0(J_2)=\{\xi:|\xi|\le1\}$,
but the only vectors $x$ occurring in generating pairs for the interior point
$\xi=0$ lie in $\C e_2$, so $m_0(J_2,0)=1$.
Indeed, if $x=(a,b)^\top$ is unit and $y\perp x$ is unit,
then $|x^*J_2y|=|a|^2$. Thus $0$ forces $a=0$.
For $n\ge3$, by contrast, the same rank-one lift and the
strict-disk criterion prove full interior multiplicity at $q=0$.
\end{remark}

\section{Boundary multiplicity for general matrices}
\label{sec:support}

Fix a boundary point and a supporting direction. We first choose a
scalar shift so that the support bound becomes an equality in the
operator norm estimate. The equality case restricts the second vector
to the maximal right singular subspace of the shifted matrix. On that
subspace the condition $x^*y=q$ becomes an ordinary numerical-range
constraint, which is needed even when the maximal singular value is
multiple.

For a compact convex set $K\subset\C$, write
\[
 h_K(\theta)=\max_{\xi\in K}\operatorname{Re}(e^{-i\theta}\xi).
\]
The corresponding supporting line is
$\operatorname{Re}(e^{-i\theta}\xi)=h_K(\theta)$.
For nonscalar $A$ and $0<q<1$, see
Li, Mehta, and Rodman~\cite[Theorem~2.8]{LiMehtaRodman} for
boundary smoothness and Chien and Nakazato~\cite[proof of Theorem~2]{ChienNakazato2002}
for strict convexity. The proof below recovers strict convexity from
the equality case of the support bound.

\subsection{The support function and its minimizing shift}
\begin{theorem}\label{thm:support-dual}
Let $A\in M_n(\C)$, $n\ge2$, and $0<q<1$.  For every
$\theta\in\mathbb R$, put $B=e^{-i\theta}A$. Then
\begin{equation}\label{eq:duality}
 h_{W_q(A)}(\theta)=
 \min_{t\in\C}\{\|B+tI\|-q\operatorname{Re}t\}.
\end{equation}
There is a unique minimizing scalar $t$. For nonscalar $A$, the
intersection of any supporting line with $W_q(A)$ is a singleton.
\end{theorem}

\begin{proof}
For an admissible pair $(x,y)$ and $t\in\C$,
\[
 \operatorname{Re}(x^*By)
 =\operatorname{Re}\bigl(x^*(B+tI)y\bigr)
       -q\operatorname{Re}t
 \le \|B+tI\|-q\operatorname{Re}t.
\]
The minimum exists since
\[
\|B+tI\|-q\operatorname{Re}t
 \ge (1-q)|t|-\|B\|\longrightarrow+\infty
 \quad(|t|\to\infty).
\]
Set $s=\sqrt{1-q^2}$, and, for a
density matrix $X\in\PSD$, write
\[
 \mu_X=\tr(BX),\qquad v_X=\tr(BB^*X),\qquad
 \delta_X=\sqrt{v_X-|\mu_X|^2}.
\]
By Cauchy--Schwarz, $\delta_X\ge0$. Put
\[
 G(X,t)=\sqrt{\tr((B+tI)(B+tI)^*X)}
                 -q\operatorname{Re}t.
\]
The function $G$ is convex in $t$ and concave in $X$, and
\[
G(X,t)=\sqrt{|\mu_X+t|^2+\delta_X^2}
                    -q\operatorname{Re}t.
\]
Put $t=-\mu_X+u+iv$, where $u,v\in\mathbb R$. Only the square root
depends on $v$, through $v^2$, so it suffices to take $v=0$.
If $\delta_X>0$, the derivative with respect to $u$ vanishes exactly when
\[
 \frac{u}{\sqrt{u^2+\delta_X^2}}=q,
 \qquad u=\frac{q\delta_X}{\sqrt{1-q^2}}.
\]
Substitution gives the following minimum:
\begin{equation}\label{eq:scalar-min}
 \inf_{t\in\C}G(X,t)
 =q\operatorname{Re}\mu_X+
   \inf_{u\in\mathbb R}(\sqrt{u^2+\delta_X^2}-qu)
 =q\operatorname{Re}\mu_X+s\delta_X.
\end{equation}
The equality remains true when $\delta_X=0$, since the infimum
of $|u|-qu$ is zero.

Sion's minimax theorem~\cite{Sion}, applied to $G$ on the compact
convex set $\PSD$ and $t\in\C$, gives
\begin{align}
\inf_{t\in\C}\max_{X\in\PSD}G(X,t)
 &=\max_{X\in\PSD}\inf_{t\in\C}G(X,t)\\
 &=\max_{X\in\PSD}
   \left(q\operatorname{Re}\mu_X+s\delta_X\right).
\label{eq:minimax}
\end{align}
Since $\max_{X\in\PSD}\tr(CC^*X)=\|C\|^2$ for $C=B+tI$,
the left side is the right side of~\eqref{eq:duality}.
By Lemma~\ref{lem:purify}, the right side of~\eqref{eq:minimax}
is
\[
\max_{\|x\|=1}
\left(q\operatorname{Re}(x^*Bx)
+s\sqrt{\|B^*x\|^2-|x^*Bx|^2}\right).
\]
The disk representation in Lemma~\ref{lem:disk-union} identifies
this maximum with $h_{W_q(A)}(\theta)$, proving~\eqref{eq:duality}.

Finally, choose a minimizing $t$ and set $\sigma=\|B+tI\|$.
For any generating pair $(x,y)$ of a point on the supporting line,
\[
\sigma=h_{W_q(A)}(\theta)+q\operatorname{Re}t
 =\operatorname{Re}(x^*(B+tI)y)
 \le |x^*(B+tI)y|\le\sigma.
\]
Equality in Cauchy--Schwarz also gives $(B+tI)y=\sigma x$.
Hence $x^*(B+tI)y=\sigma$, so
\begin{equation}\label{eq:unique-boundary-point}
 e^{-i\theta}x^*Ay=\sigma-qt.
\end{equation}
Thus each supporting line meets $W_q(A)$ in one point.

To prove uniqueness of the shift, let $t_1,t_2$ be minimizers and put
$\sigma_j=\|B+t_jI\|$. Choose one admissible pair $(x,y)$ attaining
the support maximum. Equality in both support bounds gives
\[
 (B+t_jI)y=\sigma_jx\qquad(j=1,2).
\]
Subtracting gives
\[
 (t_1-t_2)y=(\sigma_1-\sigma_2)x.
\]
Since $x,y$ are linearly independent, $t_1=t_2$.
\end{proof}

\begin{remark}[Dependence on $q$]
The support formula also shows that, for fixed $\theta$,
$q\mapsto h_{W_q(A)}(\theta)$ is concave on $(0,1)$: it is the
infimum of affine functions of $q$. Therefore
\[
 \alpha W_{q_1}(A)+(1-\alpha)W_{q_2}(A)
 \subseteq W_{\alpha q_1+(1-\alpha)q_2}(A)
 \qquad(0<q_1,q_2<1,\ 0\le\alpha\le1).
\]

\end{remark}

\subsection{The top singular subspace}
If $\|C\|=\sigma>0$ and $x,y$ are unit vectors, then
\[
 \operatorname{Re}(x^*Cy)\le |x^*Cy|
 \le\|Cy\|\le\sigma
\]
for unit vectors $x,y$. Equality throughout requires
$\|Cy\|=\sigma$, so $y$ is a right singular vector for $\sigma$,
and equality in Cauchy--Schwarz gives $x=Cy/\sigma$.
We also use the following form of Embry's multiplicity
principle~\cite{Embry}.
\begin{fixedexpectationlemma}
Let $T\in M_r(\C)$ and $w\in W(T)$. Then
\[
\dim_{\C}\operatorname{span}\{v:\|v\|=1,\ v^*Tv=w\}=r
\quad\Longleftrightarrow\quad w\in\operatorname{ri}W(T).
\]
If $w\in\operatorname{ri}W(T)$, there is also a positive definite
matrix $X$ with $\operatorname{tr}X=1$ and $\operatorname{tr}(TX)=w$.
\end{fixedexpectationlemma}
\begin{proof}
Put $\tau=\operatorname{tr}(T)/r$.
We have $\tau\in\operatorname{ri}W(T)$. Otherwise, a supporting
functional nonconstant on $W(T)$ would attain its maximum at
$\tau$. For $H=\operatorname{Re}(e^{-i\phi}T)$ this implies
\[
 \frac{\tr H}{r}
 =\operatorname{Re}(e^{-i\phi}\tau)=\lambda_{\max}(H).
\]
Hence $H$ is scalar, a contradiction. The singleton case is immediate.
For $w\in\operatorname{ri}W(T)$, choose $w_1\in W(T)$ and
$0<\varepsilon<1$ such that
$w=(1-\varepsilon)w_1+\varepsilon\tau$.
If $v_1^*Tv_1=w_1$ and $\|v_1\|=1$, then
\[
 X=(1-\varepsilon)v_1v_1^*+\varepsilon I/r>0,
 \qquad\operatorname{tr}(TX)=w.
\]
Every extreme point of the density-matrix slice
$\{Y\ge0:\operatorname{tr}Y=1,\ \operatorname{tr}(TY)=w\}$
has rank one, by the argument of Lemma~\ref{lem:purify}.
Thus $X$ is a convex combination of rank-one matrices $vv^*$
with $v^*Tv=w$. Since $X>0$, these vectors span $\C^r$.

Conversely, if $w\notin\operatorname{ri}W(T)$, choose a real
supporting functional at $w$ that is not constant on $W(T)$.
Every representing vector lies in the proper maximal eigenspace
of the associated Hermitian matrix $\operatorname{Re}(e^{-i\phi}T)$.

\end{proof}

\begin{theorem}[Exact boundary-multiplicity reduction]
\label{thm:boundary-exact}
Assume $A$ is nonscalar and $0<q<1$.  Let
$\xi_\theta$ be the unique point of $W_q(A)$ on the supporting line
with outward normal $e^{i\theta}$. Let $t=t_\theta$ be the unique
optimal shift in~\eqref{eq:duality}, and set
\[
 C_0=e^{-i\theta}A+tI,\qquad
 \sigma=\|C_0\|>0,\qquad
 E=\ker(C_0^*C_0-\sigma^2I).
\]
Let $U=\sigma^{-1}C_0|_E$ and let $P_E$ denote the
orthogonal projection onto $E$. Define
\begin{equation}\label{eq:compressed-C}
 T=P_EU^*|_E=\sigma^{-1}P_EC_0^*|_E.
\end{equation}
Then $q\in W(T)$, and
\begin{equation}\label{eq:exact-boundary}
 \xi_\theta=e^{i\theta}(\sigma-qt),\qquad
 m_q(A,\xi_\theta)
  =\dim_{\C}\operatorname{span}\{v\in E:\|v\|=1,
                         \ v^*Tv=q\}.
\end{equation}
Moreover,
\[
1\le m_q(A,\xi_\theta)\le \dim E,\]

\[
m_q(A,\xi_\theta)=\dim E \iff q\in\operatorname{ri}W(T).
\]
In particular, if $\sigma$ is simple, then
$m_q(A,\xi_\theta)=1$.
\end{theorem}

\begin{proof}
The map $U$ is an isometry, since
$\|Uv\|^2=\sigma^{-2}v^*C_0^*C_0v=\|v\|^2$ for $v\in E$.
For a generating pair $(x,y)$ at $\xi_\theta$, equality in the
support formula gives $x^*C_0y=\sigma$.  The equality case of
the singular-value inequality gives $y\in E$, $\|y\|=1$ and
$x=Uy$.  The constraint $x^*y=q$ becomes
\[
  y^*U^*y=y^*Ty=q.
\]
Conversely, each unit $y\in E$ satisfying $y^*Ty=q$ gives the
admissible pair $(Uy,y)$ and the boundary value~\eqref{eq:exact-boundary}.
Since $U$ is an isometry, it preserves the complex dimension of the
span.  This proves the formula and $q\in W(T)$.

The final equivalence follows from the fixed-expectation lemma.
\end{proof}

\subsection{A Hermitian eigenspace formula}
When $q\notin\operatorname{ri}W(T)$, the remaining condition in
Theorem~\ref{thm:boundary-exact} reduces to two Hermitian Rayleigh
quotients.

\begin{theorem}[Hermitian eigenspace formula for boundary multiplicity]
	\label{thm:boundary-hermitian}
	Let $A\in M_n(\C)$ be nonscalar and $0<q<1$.
	Fix a boundary point $\xi_\theta$ of $W_q(A)$ and let $E$, $U$, and $T$ be as in
	Theorem~\ref{thm:boundary-exact}. If $q\in\operatorname{ri}W(T)$, then
	\[
	m_q(A,\xi_\theta)=\dim E.
	\]
	Suppose that $q\notin\operatorname{ri}W(T)$.
	Choose a supporting line of $W(T)$ at $q$ with outward
	normal $e^{i\varphi}$, and set
	\[
	H_\varphi=\operatorname{Re}(e^{-i\varphi}T),\qquad
	\kappa=\lambda_{\max}(H_\varphi),\qquad
	F=\ker(H_\varphi-\kappa I_E).
	\]
	On $F$, define
	\[
	D=P_F\operatorname{Im}(e^{-i\varphi}T)|_F,
	\qquad \beta=\operatorname{Im}(e^{-i\varphi}q),
	\]
	where $P_F$ denotes the orthogonal projection onto $F$.
	Let $d_-$ and $d_+$ be the smallest and largest eigenvalues
	of the Hermitian operator $D$, respectively. Then
	$d_-\leq\beta\leq d_+$ and
	\[
	m_q(A,\xi_\theta)=
	\begin{cases}
		\dim F, & d_-<\beta<d_+,\\
		\dim\ker(D-\beta I_F),
		& \beta=d_-\text{ or }\beta=d_+.
	\end{cases}
	\]
	The second case includes $d_-=d_+=\beta$, in which case
	$m_q(A,\xi_\theta)=\dim F$.
	More precisely, the complex span of the first generating
	vectors at $\xi_\theta$ is $UF$ in the first case and
	$U\ker(D-\beta I_F)$ in the second.
\end{theorem}

\begin{proof}
	Only the case $q\notin\operatorname{ri}W(T)$ needs proof.
The supporting-line condition gives
	\[
	\operatorname{Re}(e^{-i\varphi}q)=\kappa.
	\]
	Theorem~\ref{thm:boundary-exact} guarantees that
	$q\in W(T)$, so there exists a unit vector $v\in E$ with
	$v^*Tv=q$. For any such vector,
	\[
	v^*H_\varphi v
	=\operatorname{Re}(e^{-i\varphi}v^*Tv)
	=\kappa.
	\]
	Thus $v\in F$. This also holds when $W(T)$ is a line segment
and $F=E$.
	
	For every unit vector $v\in F$, we have
	\[
	\operatorname{Re}(e^{-i\varphi}v^*Tv)
	=v^*H_\varphi v
	=\kappa
	=\operatorname{Re}(e^{-i\varphi}q),
	\]
	and
	\[
	\operatorname{Im}(e^{-i\varphi}v^*Tv)
	=v^*Dv.
	\]
	Consequently, for every unit vector $v\in F$,
	\[
	v^*Tv=q\quad\Longleftrightarrow\quad v^*Dv=\beta.
	\]
	Since $q\in W(T)$, the Hermitian Rayleigh quotient gives
	\[
	d_-\leq\beta\leq d_+.
	\]
	If $\beta=d_+$ or $d_-$, the generating vectors span the
corresponding extreme eigenspace of $D$. This also covers
$d_-=d_+=\beta$.
	
	Now suppose $d_-<\beta<d_+$. Let
	$f_1,\ldots,f_s$ be an orthonormal eigenbasis of $D$,
	with corresponding eigenvalues $d_1,\ldots,d_s$.
	There exist strictly positive numbers $p_1,\ldots,p_s$
	such that
	\[
	\sum_{j=1}^s p_j=1,\qquad
	\sum_{j=1}^s p_jd_j=\beta.
	\]
	Such weights exist because $d_-<\beta<d_+$.
	For every $\varepsilon=(\varepsilon_1,\ldots,\varepsilon_s)
	\in\{-1,1\}^s$, put
	\[
	v_\varepsilon
	=\sum_{j=1}^s\varepsilon_j\sqrt{p_j}\,f_j.
	\]
	Each $v_\varepsilon$ is a unit vector with
$v_\varepsilon^*Dv_\varepsilon=\beta$, and
	\[
	f_j=\frac{1}{2^s\sqrt{p_j}}
	\sum_{\varepsilon\in\{-1,1\}^s}
	\varepsilon_j v_\varepsilon.
	\]
	Therefore the $v_\varepsilon$ span $F$. The assertion now
follows from Theorem~\ref{thm:boundary-exact} by applying $U$.
\end{proof}

\begin{example}[A second eigenspace calculation may be needed]
\label{ex:boundary-two-step}
Let $q=\tfrac12$ and consider the $4\times4$ matrix
\[
 A=\begin{pmatrix}
 \frac12&0&0&0\\
 0&\frac12-\frac{i}{4}&0&0\\
 \frac{\sqrt3}{2}&0&0&0\\
 0&\frac{\sqrt{11}}4&0&0
 \end{pmatrix}.
\]
The first two columns of $A$ are orthonormal and the other
columns vanish. Hence $\|A\|=1$, and its top right singular
space is $E=\spn\{e_1,e_2\}$. The vectors $y=e_1$ and
$x=Ae_1$ satisfy $x^*y=\tfrac12$ and $x^*Ay=1$.
Consequently $1\in\partial W_{1/2}(A)$, since
$|x^*Ay|\le\|A\|=1$ for every admissible pair.
The shift $t=0$ is optimal in the outward direction $\theta=0$.
The compression of Theorem~\ref{thm:boundary-exact} is
\[
             T=P_EA^*|_E
              =\operatorname{diag}\!\left(\frac12,
                                \frac12+\frac{i}{4}\right).
\]
Its numerical range is the vertical segment from $\frac12$ to
$\frac12+\frac{i}{4}$. Choosing the supporting line
$\operatorname{Re}w=\frac12$ gives $F=E$, so the first Hermitian
eigenspace alone does not determine the multiplicity. However,
\[
              D=\operatorname{Im}T
                 =\operatorname{diag}(0,\tfrac14),
              \qquad\beta=0.
\]
Theorem~\ref{thm:boundary-hermitian} now yields
\[
               m_{1/2}(A,1)=\dim\ker D=1,
\]
although $\dim E=2$. 
\end{example}

\begin{remark}
	The optimal shift is unique by Theorem~\ref{thm:support-dual}.
	Consequently, the compressed operator $T$ is determined up to
	unitary equivalence by $A$, $q$, and the supporting direction.
\end{remark}

\begin{corollary}
\label{cor:normality-obstruction}
Let $A\in M_n(\C)$ be nonscalar, $0<q<1$, and
$\xi\in\partial W_q(A)$. If $m_q(A,\xi)=n$, then $A$ is normal.
More precisely, for the optimal shift in Theorem~\ref{thm:boundary-exact},
\[
             e^{-i\theta}A+tI=\sigma V
\]
with $V$ unitary and $q\in\operatorname{ri}W(V^*)$.
Conversely, if the optimal shift has this property, then
$m_q(A,\xi)=n$. In particular, if $A$ is nonnormal,
\[
           m_q(A,\xi)\le n-1\qquad(\xi\in\partial W_q(A)).
\]
\end{corollary}
\begin{proof}
Theorem~\ref{thm:boundary-exact} gives
$m_q(A,\xi)\le \dim E\le n$. If $m_q(A,\xi)=n$, the top singular
subspace $E$ is all of $\C^n$, so
$C_0^*C_0=\sigma^2 I$ and $C_0=\sigma V$ for a unitary $V$.
It follows that $A$ is normal. The criterion
$q\in\operatorname{ri}W(V^*)$ follows from the equality statement
of that theorem. The converse is the same equality statement.
\end{proof}

\section{Hermitian and normal matrices}\label{sec:hermitian}
For Hermitian matrices, $W_q(A)$ is an ellipse~\cite{LiNakazato1998,Tsing}.
We determine the multiplicity in terms of the extreme eigenspaces.
\begin{theorem}[Hermitian boundary multiplicity]\label{thm:herm}
Let $H=H^*$ be nonscalar, let $a<b$ be its extreme eigenvalues,
and write $c=(a+b)/2$, $\Delta=(b-a)/2$. Then
\[
W_q(H)=\left\{\xi:\frac{(\operatorname{Re}\xi-qc)^2}{\Delta^2}
                +\frac{(\operatorname{Im}\xi)^2}{s^2\Delta^2}\le1\right\},
\]
and for every $\xi\in\partial W_q(H)$,
\[
             m_q(H,\xi)=\dim E_a+\dim E_b,
\]
where $E_a,E_b$ are the eigenspaces of $H$ for $a,b$.
\end{theorem}
\begin{proof}
By covariance reduce to $K=(H-cI)/\Delta$, a Hermitian contraction
with eigenvalues $-1$ and $1$. For an admissible pair write
$y=qx+su$ with $x^*u=0$ and put
$\alpha=x^*Kx\in[-1,1]$, $\eta=x^*Ku$. Since $K^2\le I$,
\[
 |\eta|^2\le\|Kx\|^2-|x^*Kx|^2\le1-\alpha^2.
\]
Therefore $x^*Ky=q\alpha+s\eta$ lies in the union of disks
$|w-q\alpha|\le s\sqrt{1-\alpha^2}$, $-1\le\alpha\le1$.
That union is the ellipse
$ (\operatorname{Re}w)^2+(\operatorname{Im}w)^2/s^2\le1$:
its support function in direction $(\cos\theta,\sin\theta)$ is
\[
\max_{-1\le\alpha\le1}
\left(q\alpha\cos\theta+s\sqrt{1-\alpha^2}\right)
=\sqrt{q^2\cos^2\theta+s^2}.
\]
Since $q^2+s^2=1$,
\[
 q^2\cos^2\theta+s^2=\cos^2\theta+s^2\sin^2\theta.
\]
Thus this support function equals that of the stated ellipse,
which proves that the union of disks is contained in the ellipse.
The reverse inclusion follows from the explicit construction below.
Conversely, choose any unit $e_+\in E_1$ and $e_-\in E_{-1}$.
For a prescribed $\alpha\in[-1,1]$, put
\[
 x=\sqrt{\frac{1+\alpha}{2}}e_+
   +\sqrt{\frac{1-\alpha}{2}}e_-,\quad
 u=e^{i\phi}\left(\sqrt{\frac{1-\alpha}{2}}e_+
    -\sqrt{\frac{1+\alpha}{2}}e_-\right).
\]
Then $x^*Kx=\alpha$ and
$x^*Ku=e^{i\phi}\sqrt{1-\alpha^2}$.
We verify that these circles fill the ellipse.  Write
$w=X+iY$.  The circle condition is
\[
  (X-q\alpha)^2+Y^2=s^2(1-\alpha^2),
\]
or, since $q^2+s^2=1$,
\[
  (\alpha-qX)^2=s^2(1-X^2)-Y^2.
\]
If $X^2+Y^2/s^2\le1$, the right side is nonnegative, and
we may choose
$\alpha=qX+\sqrt{s^2(1-X^2)-Y^2}$.  Moreover,
\[
 |\alpha|\le q|X|+s\sqrt{1-X^2}\le1
\]
by Cauchy--Schwarz.  Thus a permitted $\alpha\in[-1,1]$
exists, and an appropriate phase $\phi$ gives $w$.
The reverse inclusion follows from the same circle equation.

At a boundary value $w$, choose its supporting direction
$e^{i\theta}$. For every generating pair,
\begin{align*}
\operatorname{Re}(e^{-i\theta}w)
&\le q\alpha\cos\theta+s\sqrt{\|Kx\|^2-\alpha^2}\\
&\le q\alpha\cos\theta+s\sqrt{1-\alpha^2}
\le h_{W_q(K)}(\theta).
\end{align*}
The first and last expressions are equal at $w$.
Since $s>0$, equality in the middle inequality forces
$\|Kx\|=1$. Hence
\[
 0=\|x\|^2-\|Kx\|^2=x^*(I-K^2)x.
\]
Since $I-K^2$ is positive semidefinite, this forces
$(I-K^2)x=0$ and $x\in E_1\oplus E_{-1}$. Thus boundary multiplicity is at most
$\dim E_1+\dim E_{-1}$.
The ellipse has a unique supporting direction at $w$; maximizing
$q\alpha\cos\theta+s\sqrt{1-\alpha^2}$ gives
$\alpha=q\cos\theta/\sqrt{q^2\cos^2\theta+s^2}
=q\operatorname{Re}w$, and yields
$|\alpha|<1$ (as $0<q<1$). Hence the two components of the
constructed $x$ are both nonzero. Varying $e_+$ and $e_-$
and taking linear combinations of such generating vectors spans
all of $E_1\oplus E_{-1}$. This proves the reverse inequality.
\end{proof}
\begin{corollary}[Affine Hermitian matrices]\label{cor:affine-hermitian}
If $A=\gamma H+\delta I$ with $H=H^*$ nonscalar and
$\gamma\ne0$, every interior point has multiplicity $n$, and
every boundary point has multiplicity $\dim E_a+\dim E_b$.
\end{corollary}
\begin{proof}Apply Proposition~\ref{prop:inv},
Theorems~\ref{thm:all} and~\ref{thm:herm}.\end{proof}
\begin{remark}[Normal matrices]
For $A=\operatorname{diag}(\lambda_1,\ldots,\lambda_n)$, writing
$p_j=|x_j|^2$ reduces the strict-disk condition to
\[
s^2\sum_jp_j|\lambda_j|^2
 -\Big|\sum_jp_j\lambda_j-q\xi\Big|^2-s^2|\xi|^2>0.
\]
These weights encode the spectral geometry studied by Chien and
Nakazato~\cite{ChienNakazato2005}.
The general lifting proof replaces the probability vector $(p_j)$
by a density matrix $X$, and Lemma~\ref{lem:purify} removes the
normality assumption. See also~\cite{Li1998,LiNakazato1998} for earlier results
on $q$-numerical ranges.
\end{remark}

The general normality obstruction, together with the Hermitian
formula, gives a particularly simple result in dimension two.

\begin{corollary}
\label{cor:two-by-two-positive}
Let $A\in M_2(\C)$ be nonscalar and $0<q<1$. Then
\[
 m_q(A,\xi)=
 \begin{cases}
 2,&\xi\in\operatorname{int}W_q(A),\\
 2,&\xi\in\partial W_q(A)\text{ and }A\text{ is normal},\\
 1,&\xi\in\partial W_q(A)\text{ and }A\text{ is nonnormal}.
 \end{cases}
\]
In particular, every boundary point has multiplicity two if and
only if $A$ is normal.
\end{corollary}
\begin{proof}
The interior assertion is Theorem~\ref{thm:all}.
A nonscalar normal $2\times2$ matrix has two distinct eigenvalues,
and it is an affine transformation of a Hermitian matrix with
extreme eigenspaces of dimension one. Its boundary multiplicity is
two by Corollary~\ref{cor:affine-hermitian}. If $A$ is nonnormal,
Corollary~\ref{cor:normality-obstruction} gives
$m_q(A,\xi)\le1$ at any boundary point, while every point of
$W_q(A)$ has at least one generating pair.
\end{proof}

For normal matrices of larger order, the optimal shift selects the
eigenvalues of maximal modulus. The remaining condition is a convex
combination of these active eigenvalues. A nonempty convex subset $F$
of a convex set $K$ is a \emph{face} if, whenever $a,b\in K$,
$0<t<1$, and $ta+(1-t)b\in F$, both $a$ and $b$ belong to $F$.
The whole set is allowed as an improper face. A face obtained by
intersection with a supporting line is called \emph{exposed}.
All nonempty proper faces of the finite polygons used below are
exposed, but general compressed numerical ranges in Section~\ref{sec:proj}
may also have nonexposed faces. For a point of a finite-dimensional
convex set there is a unique minimal face containing it, and the point
belongs to the relative interior of that face. The smallest face
containing $q$ identifies precisely which eigenspaces can occur in a
generating pair.

\begin{theorem}
\label{thm:normal-face}
Let $A$ be nonscalar and normal, with eigenvalues
$\lambda_1,\ldots,\lambda_n$ counted with multiplicity.
Fix $0<q<1$, $\theta\in\mathbb R$, and let $t$ be the minimizer
in~\eqref{eq:duality}. Put
\[
 \sigma=\max_j|e^{-i\theta}\lambda_j+t|,
 \qquad S=\{j:|e^{-i\theta}\lambda_j+t|=\sigma\},
\]
\[
 \omega_j=(e^{-i\theta}\lambda_j+t)/\sigma\quad(j\in S),
 \qquad P=\operatorname{conv}\{\omega_j:j\in S\}.
\]
Let $F$ be the unique smallest face of $P$ with
$q\in\operatorname{ri}F$, viewed in $\mathbb R\subset\C$.
Then
\[
m_q(A,\xi_\theta)=\#\{j\in S:\omega_j\in F\},
\]
where $\xi_\theta$ is the unique supported point. The equality
$m_q(A,\xi_\theta)=|S|$ holds if and only if $q\in\operatorname{ri}P$.
\end{theorem}
\begin{proof}
In a unitary eigenbasis of $A$, the subspace $E$ of
Theorem~\ref{thm:boundary-exact} is $\spn\{e_j:j\in S\}$.
On this basis $Ue_j=\omega_j e_j$ and $Te_j=\overline{\omega_j}e_j$.
Since $q$ is real, the conditions $\|v\|=1$ and $v^*Tv=q$ reduce to
\[
 \sum_{j\in S}|v_j|^2=1,\qquad
 \sum_{j\in S}|v_j|^2\omega_j=q.
\]
The latter identity expresses $q$ as a convex combination in $P$.
The face property gives $v_j=0$ whenever $\omega_j\notin F$.
For the reverse inclusion, let $J=\{j\in S:\omega_j\in F\}$ and
let $w=|J|^{-1}\sum_{j\in J}\omega_j$. Since $q\in\operatorname{ri}F$,
we have $(q-\varepsilon w)/(1-\varepsilon)\in F$ for sufficiently
small $\varepsilon>0$. Expressing this point as a convex combination
of the $\omega_j$ with $j\in J$ yields a representation of $q$ with
coefficients $p_j\ge\varepsilon/|J|$ on $J$. Choose $v_j=\sqrt{p_j}$
on $J$ and zero elsewhere. Changing the sign of one coordinate at a
time preserves both constraints and, by subtraction, puts every $e_j$
with $j\in J$ in the span. Theorem~\ref{thm:boundary-exact} now gives
the asserted dimension.
\end{proof}

\begin{corollary}[A disk range with full boundary multiplicity]
\label{cor:unitary-polygon}
Let $V\in M_n(\C)$ be unitary, and suppose its eigenvalue polygon
$P=\operatorname{conv}\sigma(V)$ contains the closed disk
$\{w:|w|\le\rho\}$ for some $0<\rho<1$.
For every $0<q\le\rho$,
\[
                         W_q(V)=\{\xi:|\xi|\le1\}.
\]
If $0<q<\rho$, then $m_q(V,\xi)=n$ for every $|\xi|=1$.
For $q=\rho$, use the face formula in Theorem~\ref{thm:normal-face}
with $t=0$.
\end{corollary}
\begin{proof}
As $V$ is unitary, $|x^*Vy|\le1$. The inclusion hypothesis gives
$qe^{i\theta}\in W(V)=P$ for every $\theta$. Choose a unit $y$ with
$y^*Vy=qe^{i\theta}$ and set $x=e^{-i\theta}Vy$. Directly,
$x^*y=q$ and $x^*Vy=e^{i\theta}$. Thus the whole unit circle belongs
to $W_q(V)$, and convexity gives the closed disk.
The support value is $1$, attained by the shift $t=0$ in
\eqref{eq:duality}. All eigenvalues are therefore active. If $q<\rho$,
$q e^{i\theta}$ is interior to $P$, so Theorem~\ref{thm:normal-face}
gives $m_q(V,e^{i\theta})=n$. For $q=\rho$ the same face formula
applies whether the point is interior or on an edge of $P$.
\end{proof}
\begin{example}[Three cube roots of unity]
\label{ex:roots-unity}
Let $\omega=e^{2\pi i/3}$ and
$A=\operatorname{diag}(1,\omega,\omega^2)$.
The triangle $\operatorname{conv}\{1,\omega,\omega^2\}$ has inradius
$1/2$, so Corollary~\ref{cor:unitary-polygon} gives the unit disk for
$0<q\le1/2$. If $q<1/2$, all three indices contribute at every
boundary point. When $q=1/2$, the circle of radius $q$ touches the
three edges at their midpoints. At those three supporting directions
$m_q=2$; at every other direction $m_q=3$. Here the transition in
multiplicity occurs without a change in the range.
\end{example}

\begin{example}[Full multiplicity can occur on the boundary]
For $H=\operatorname{diag}(-1,1)$, Theorem~\ref{thm:herm}
gives $m_q(H,\xi)=2$ throughout the ellipse, including its boundary.
Inserting the additional eigenvalue $0$ gives
$H=\operatorname{diag}(-1,0,1)$, with the same ellipse and boundary
multiplicity $2$, but interior multiplicity $3$. In neither example
is strict convexity accompanied by a simple maximal singular value
at the boundary.
\end{example}

\section{Weighted nilpotent shifts}\label{sec:weighted}
The $q$-numerical radii of positive weighted nilpotent shifts were
studied in~\cite[Proposition~4.3]{LiNakazato1998}
and~\cite{ChienNakazato2002,ChienNakazato,ChienNakazatoPerturbation}.
Here the positions of the maximal weights determine the multiplicities.

Let
\[
A=\begin{pmatrix}
0&a_1&&&\\&0&a_2&&\\&&\ddots&\ddots&\\&&&0&a_{n-1}\\&&&&0
\end{pmatrix},\qquad a_j>0,
\]
let $M=\max_j a_j$, and let $S=\{j:a_j=M\}$. Decompose $S$ into maximal consecutive runs of lengths $\ell_1,\ldots,\ell_k$, and put
\[
r=|S|,\qquad L=\max_i\ell_i,\qquad k_L=\#\{i:\ell_i=L\},\qquad q_A=\cos\frac{\pi}{L+1}.
\]
A run of $\ell$ maximal weights contributes $\ell$ coordinate
directions to the maximal right singular subspace; the Hermitian part
of the associated compressed shift has top eigenvalue
$\cos(\pi/(\ell+1))$. Put
$D_\theta=\operatorname{diag}(e^{i(j-1)\theta})_{j=1}^n$.
Since $D_\theta^*AD_\theta=e^{i\theta}A$, covariance and convexity
imply that $W_q(A)$ is a disk. Theorem~\ref{thm:all} gives
$m_q(A,\xi)=n$ in its interior.

\begin{theorem}[Boundary multiplicity of weighted shifts]\label{thm:weighted}
For the above positive weighted nilpotent shift and every $\xi\in\partial W_q(A)$,
\[
m_q(A,\xi)=\begin{cases}
r,&0<q<q_A,\\
k_L,&q=q_A>0,\\
1,&q_A<q<1.
\end{cases}
\]
Moreover, writing $f(t)=\|A+tI\|$, the radius satisfies
\[
R_q=\begin{cases}M,&0<q\le q_A,\\ f(t_q)-qt_q,&q_A<q<1,\end{cases}
\]
where $t_q>0$ is the unique solution of $f'(t_q)=q$.
If $q_A=0$, the first two multiplicity cases are empty for
$0<q<1$, and only $m_q(A,\xi)=1$ occurs on the boundary.
\end{theorem}
\begin{proof}
Because $|x^*Ay|\le\|Ay\|\le M$, equality $x^*Ay=M$ holds precisely when $y$ is a unit vector in
\[
E=\spn\{e_{j+1}:j\in S\},\qquad x=Ay/M,
\]
and $(Ay)^*y/M=q$. The restriction $A/M:E\to AE$ is an isometry.
In the coordinate basis of $E$, two consecutive selected indices
are coupled by an entry $1$ of $P_EA|_E/M$; gaps separate the blocks.
Taking the Hermitian part gives
\[
H_E=T_{\ell_1}\oplus\cdots\oplus T_{\ell_k},\qquad
T_\ell=\frac12\begin{pmatrix}0&1&&\\1&0&\ddots&\\&\ddots&\ddots&1\\&&1&0\end{pmatrix}.
\]
Its extreme eigenvalues are $\pm q_A$ and its top eigenspace has dimension $k_L$. The real Rayleigh quotients of $H_E$ fill $[-q_A,q_A]$, and for real vectors the skew-symmetric part of the compression contributes zero.
Consequently, for each $0<q\le q_A$ there is a real unit vector
$y\in E$ with
\[
 (Ay)^*y/M=y^*H_Ey=q.
\]
Taking $x=Ay/M$ gives $\|x\|=\|y\|=1$,
$x^*y=q$, and $x^*Ay=M$. As $|x^*Ay|\le\|A\|=M$ always,
the boundary radius is $R_q=M$ in this interval.

For $0<q<q_A$, take a real orthonormal eigenbasis $u_1,\ldots,u_r$ of $H_E$ with eigenvalues $\lambda_j$. As $q$ lies strictly between the extreme eigenvalues, choose strictly positive $p_j$ with $\sum p_j=1$ and $\sum p_j\lambda_j=q$. All sign choices
\[
y_\varepsilon=\sum_{j=1}^r\varepsilon_j\sqrt{p_j}\,u_j,\qquad \varepsilon_j\in\{-1,1\},
\]
are admissible second vectors generating $M$, and their span is $E$. The corresponding $x$-vectors span $AE$, also of dimension $r$.

At $q=q_A>0$, equality in the top Rayleigh quotient forces $y$ into the top eigenspace of $H_E$. On each longest run its top eigenvector is real, and the different runs are disconnected in the compression, so the skew-Hermitian part has zero quadratic form on their entire span. Hence all unit vectors in that $k_L$-dimensional eigenspace are admissible, giving multiplicity $k_L$.

For $t>0$, $(A+tI)^*(A+tI)$ is nonnegative irreducible tridiagonal, with strictly positive off-diagonal entries $ta_j$. Perron--Frobenius theory gives a simple top eigenvalue for
this irreducible nonnegative Hermitian matrix.  Thus $f(t)$ is
real analytic on $(0,\infty)$.  Choose the corresponding
left and right top singular vectors $u(t),v(t)$ to be real,
strictly positive, and of unit length.  The simple-singular-value
perturbation formula gives
\[
f'(t)=\operatorname{Re}(u(t)^*v(t))=u(t)^*v(t).
\]
Since $(A+tI)^*(A+tI)=A^*A+t(A+A^*)+t^2I$,
writing $\lambda(t)=f(t)^2$, we have $\lambda(0)=M^2$
and $\lambda_+'(0)=2M f_+'(0)$. The first-order Hermitian
eigenvalue perturbation formula at the possibly multiple top
eigenvalue of $A^*A$ therefore gives
\[
 f'_+(0)=\frac{1}{2M}\lambda_{\max}
      \big(P_E(A+A^*)|_E\big)
       =\lambda_{\max}(H_E)=q_A.
\]
Here we use the standard one-sided eigenvalue perturbation
formula~\cite[Chapter~II]{Kato}: the derivative of the top
eigenvalue at a multiple eigenvalue is the largest eigenvalue
of the perturbation compressed to its top eigenspace.
The term $t^2I$ has zero derivative at the origin.

Since $f(t)/t=\|I+A/t\|\to1$, convexity gives
\[
 \frac{f(t)-f(0)}{t}\le f'(t)
 \le\frac{f(2t)-f(t)}{t}
\]
show that $f'(t)\to1$. Convexity makes $f'$ nondecreasing.
Equality of $f'$ at two distinct positive arguments would make
$f$ affine there and hence, by analyticity, on $(0,\infty)$,
contrary to $f'_+(0)=q_A<1$ and $f'(t)\to1$.
Thus $f'$ is strictly increasing from $q_A$ to $1$.

For any admissible pair,
\[
\operatorname{Re}(x^*Ay)\le\|A+t_qI\|-t_qq.
\]
By construction the unit singular pair $(u(t_q),v(t_q))$
satisfies
\[
 u(t_q)^*v(t_q)=f'(t_q)=q.
\]
It is therefore an admissible generating pair, and equality
holds in the support bound.  Hence $R_q=f(t_q)-t_qq$.
Every pair attaining $R_q$ is a top singular pair of
$A+t_qI$, unique up to a common phase.  The span of its
generating vectors is one-dimensional.  Rotational similarity
transfers the result from the positive boundary point to the
entire circle.
\end{proof}
For the Jordan block $J_n$, the radius threshold
$q=\cos(\pi/n)$ is known~\cite[Corollary~7]{ChienNakazato2002};
see also~\cite[Proposition~4.3]{LiNakazato1998} for $J_3$.
Theorem~\ref{thm:weighted} gives the following multiplicities.
\begin{corollary}\label{cor:jordan}
For $n\ge 2$ and $0<q<1$, let $R_q=\max\{|\xi|:\xi\in W_q(J_n)\}$. Then
\[
m_q(J_n,\xi)=\begin{cases}
n,&|\xi|<R_q,\\
n-1,&|\xi|=R_q,\ 0<q<\cos(\pi/n),\\
1,&|\xi|=R_q,\ \cos(\pi/n)\le q<1.
\end{cases}
\]
\end{corollary}

\begin{proof}
Take $a_j=1$ for every $j$. Then $r=L=n-1$, $k_L=1$, and $q_A=\cos(\pi/n)$. The interior case follows from Theorem~\ref{thm:all}, and the boundary cases from Theorem~\ref{thm:weighted}.
\end{proof}

\begin{example}\label{ex:runs}
Let $A$ have weights $(2,2,1,2,2)$. Then $S=\{1,2,4,5\}$,
$r=4$, $L=2$, $k_L=2$, and $q_A=1/2$. By
Theorem~\ref{thm:weighted}, for every $\xi\in\partial W_q(A)$,
\[
m_q(A,\xi)=
\begin{cases}
4,&0<q<1/2,\\
2,&q=1/2,\\
1,&1/2<q<1.
\end{cases}
\]
\end{example}

\begin{example}[A unique maximal weight]\label{ex:unique}
Consider the $4\times4$ weighted shift with weights $(1,2,1)$. Its only maximal weight occurs at index $2$. Thus $r=L=k_L=1$ and $q_A=0$. Theorem~\ref{thm:weighted} shows that $m_q(A,\xi)=1$ at every boundary point and for every $0<q<1$. In contrast with $J_4$, there is no positive transition parameter.
\end{example}

\begin{example}[One long run and one short run]\label{ex:unequal}
Let $A$ have weights $(3,3,3,1,3,3,1)$. The maximal-weight
runs have lengths $3$ and $2$, so $r=5$, $L=3$, $k_L=1$, and
$q_A=1/\sqrt2$. By Theorem~\ref{thm:weighted}, for every
$\xi\in\partial W_q(A)$,
\[
m_q(A,\xi)=
\begin{cases}
5,&0<q<1/\sqrt2,\\
1,&1/\sqrt2\le q<1.
\end{cases}
\]
The shorter run contributes to the subcritical multiplicity but not to the multiplicity at the threshold.
\end{example}

\begin{remark}[Combinatorial rigidity and metric variation]
Changing the nonmaximal weights, while keeping them positive and
strictly below the maximum, leaves $S$, $r$, $L$, $k_L$, and $q_A$
unchanged. Hence it leaves the boundary multiplicity unchanged.
For $q>q_A$, however, the radius involves $\|A+t_qI\|$ and in general
depends on the nonmaximal weights as well.
\end{remark}

\section{Multiplicities of generating projections}\label{sec:proj}
We now consider the real span of the generating projections $xx^*$.
Unlike vector multiplicity, projection multiplicity can distinguish
between two fibers whose generating vectors span the same complex
subspace. For interior points in dimension at least three, we show
that the projections span all Hermitian matrices. At boundary points,
the singular-space reduction imposes a fixed numerical expectation
on their rank-one projections, leading to a correction to the square
of the vector multiplicity. That correction is governed by the
smallest face of a compressed numerical range.

\begin{theorem}
\label{thm:interior-projection}
Let $n\ge3$, $A\in M_n(\C)$, and $0<q<1$. For every
$\xi\in\operatorname{ri}W_q(A)$ one has
\[
                 p_q(A,\xi)=n^2.
\]
The same conclusion holds at interior points for $q=0$ and $n\ge3$.
For a scalar matrix the conclusion holds also when $n=2$ and
$0\le q<1$.
\end{theorem}
\begin{proof}
Suppose first that $A$ is nonscalar and $0<q<1$.
By Theorem~\ref{thm:strict-all} and continuity, a nonempty open
subset $\mathcal O$ of the unit sphere satisfies the strict-disk
inequality for $\xi$. Since $n\ge3$, $\mathcal O\subseteq\mathcal M_q(A,\xi)$.
Let $H=H^*$ be in the trace-orthogonal complement of the real span
of the generating projections. Then
\[
 \tr(Hxx^*)=x^*Hx=0\qquad(x\in\mathcal O).
\]
By homogeneity the quadratic form vanishes on the open cone
$\{tx:t>0,\ x\in\mathcal O\}$ in $\mathbb R^{2n}$. A real polynomial
vanishing on a nonempty open set is identically zero. Polarization
therefore gives $H=0$, proving the dimension assertion.
For $q=0$, the strict-radius condition in the proof of
Theorem~\ref{thm:q0interior} gives the required open subset.
If $A$ is scalar, every unit vector generates the sole point of
$W_q(A)$, and the same trace-orthogonality argument applies.
\end{proof}
\begin{remark}
The open-set argument requires $n\ge3$, since for $n=2$ the
fixed-first-vector values form a circle.
\end{remark}

\begin{theorem}
\label{thm:boundary-projection}
Let $A$ be nonscalar, $0<q<1$, and let $\xi_\theta$ be the
boundary point and $E,T,U$ the objects in
Theorem~\ref{thm:boundary-exact}. Define
\[
 \mathcal V=\operatorname{span}_{\C}\{v\in E:\|v\|=1,\ v^*Tv=q\},
 \quad m=\dim\mathcal V=m_q(A,\xi_\theta),
\]
let $T_{\mathcal V}=P_{\mathcal V}T|_{\mathcal V}$, and put
\[
 r_{\mathcal V}=\dim_{\mathbb R}\operatorname{span}_{\mathbb R}
 \{I_{\mathcal V},\operatorname{Re}T_{\mathcal V},\operatorname{Im}T_{\mathcal V}\}.
\]
Then
\[
p_q(A,\xi_\theta)=m^2-r_{\mathcal V}+1,
\qquad 1\le r_{\mathcal V}\le3.
\]
In particular, $m^2-2\le p_q(A,\xi_\theta)\le m^2$.
\end{theorem}
\begin{proof}
The generating first vectors are $Uv$, where $v\in\mathcal V$
is unit and $v^*T_{\mathcal V}v=q$.
The fixed-expectation lemma gives a positive definite density
matrix $X_0$ on $\mathcal V$ with
$\operatorname{tr}(T_{\mathcal V}X_0)=q$.
The affine space
\[
 \mathcal L_q=\{X=X^*\in M_m(\C):
       \operatorname{tr}X=1,\ \operatorname{tr}(T_{\mathcal V}X)=q\}
\]
has real dimension $m^2-r_{\mathcal V}$. The real constraints are
\[
 \operatorname{tr}X=1,\qquad
 \operatorname{tr}((\operatorname{Re}T_{\mathcal V})X)=q,\qquad
 \operatorname{tr}((\operatorname{Im}T_{\mathcal V})X)=0.
\]
The coefficient matrices have real span of dimension $r_{\mathcal V}$.
Since $X_0>0$, the positive semidefinite part of $\mathcal L_q$
has full affine hull $\mathcal L_q$. Its extreme points have
rank one by Lemma~\ref{lem:purify}. Thus the rank-one projections
span a real space of dimension $m^2-r_{\mathcal V}+1$.
Conjugation by the isometry $U$ preserves this dimension.
\end{proof}

The preceding formula has a geometric interpretation. The
projections supported on $\mathcal V$ have trace one and satisfy
$\operatorname{tr}(T_{\mathcal V}X)=q$. Although the latter is a
complex equation, its real and imaginary parts need not impose two
independent constraints: their number depends on the affine geometry
of the numerical range of the compression. The next result expresses
that dependence through a minimal face. Recall that a face $F$ of a
convex set $K$ has the property that, if an interior point of a segment
in $K$ belongs to $F$, then both endpoints belong to $F$.
Every point of a finite-dimensional convex set belongs to the
relative interior of a unique smallest face; the face here need
not be exposed. We next express $r_{\mathcal V}$ in terms of that
face of $W(T)$.

\begin{theorem}[Boundary projection multiplicity and faces]
\label{thm:boundary-face-projection}
Let $A$ be nonscalar, $0<q<1$, and
$\xi_\theta\in\partial W_q(A)$. Let $T$ be as in
Theorem~\ref{thm:boundary-exact}, and let $F_q$ be the smallest
face of $W(T)$ with $q\in\operatorname{ri}F_q$. Put
$m=m_q(A,\xi_\theta)$ and
$d=\dim_{\mathbb R}\operatorname{aff}F_q$. Then
\[
p_q(A,\xi_\theta)=m^2-d.
\]
\end{theorem}
\begin{proof}
Let
\[
 S=\operatorname{span}_{\C}\{v\in E:\|v\|=1,\ v^*Tv=q\},
 \qquad T_S=P_ST|_S.
\]
We claim that
\begin{equation}\label{eq:face-compression}
	W(T_S)=F_q.
\end{equation}

To prove $F_q\subseteq W(T_S)$, fix $w\in F_q$.
Because $q\in\operatorname{ri}F_q$, there exist $w'\in F_q$
and $0<\alpha<1$ with
$q=\alpha w+(1-\alpha)w'$. Choose unit vectors $v,u\in E$
satisfying $v^*Tv=w$ and $u^*Tu=w'$. The density matrix
\[
                  X=\alpha vv^*+(1-\alpha)uu^*
\]
satisfies $\operatorname{tr}(TX)=q$.
The rank-one extreme-point argument of Lemma~\ref{lem:purify}
shows that every density matrix with this fixed expectation is
a convex combination of rank-one matrices $aa^*$ with $a^*Ta=q$.
All such vectors $a$ belong to $S$. Consequently
$\operatorname{ran}X\subseteq S$, and hence $v\in S$.
Thus $w=v^*T_Sv\in W(T_S)$.

Conversely, select unit vectors
$v_1,\ldots,v_m$ with $v_j^*Tv_j=q$ which form a basis of $S$.
The density matrix
\[
                     X_0=\frac1m\sum_{j=1}^m v_jv_j^*
\]
is positive definite on $S$ and satisfies
$\operatorname{tr}(TX_0)=q$.
If $v\in S$ is any unit vector, choose $\varepsilon>0$
small enough that $X_0-\varepsilon vv^*\ge0$ and
$\varepsilon<1$. Then
\[
 X_0=\varepsilon vv^*+(1-\varepsilon)Y,
 \qquad Y\ge0,\quad\operatorname{tr}Y=1.
\]
It follows that
\[
 q=\varepsilon v^*Tv+(1-\varepsilon)\operatorname{tr}(TY).
\]
Both terms on the right belong to $W(T)$, since numerical
ranges contain expectations of density matrices. As $F_q$ is
a face containing $q$, both terms belong to $F_q$.
Therefore $W(T_S)\subseteq F_q$, proving~\eqref{eq:face-compression}.

Finally,
\[
\dim_{\mathbb R}\operatorname{span}_{\mathbb R}
 \{I_S,\operatorname{Re}T_S,\operatorname{Im}T_S\}=1+d.
\]
Indeed, the expectation map on trace-one Hermitian matrices
has affine image $\operatorname{aff}W(T_S)=\operatorname{aff}F_q$
of dimension $d$.
Theorem~\ref{thm:boundary-projection} now gives
$p_q(A,\xi_\theta)=m^2-(1+d)+1=m^2-d$.
\end{proof}

\begin{remark}
In Theorem~\ref{thm:boundary-face-projection}, $d\in\{0,1,2\}$.
If $q\notin\operatorname{ri}W(T)$, the notation of
Theorem~\ref{thm:boundary-hermitian} gives $d=1$ when
$d_-<\beta<d_+$, and $d=0$ at an endpoint.
Otherwise $d=\dim_{\mathbb R}\operatorname{aff}W(T)$.
\end{remark}

\begin{corollary}\label{cor:normal-projections}
Let $A$ be a nonscalar normal matrix, and use the notation
$S$, $P$, and $F$ of Theorem~\ref{thm:normal-face}. Put
$m=\#\{j\in S:\omega_j\in F\}$ and
$d=\dim_{\mathbb R}\operatorname{aff}F$. For $0<q<1$,
\[
        p_q(A,\xi_\theta)=m^2-d.
\]
Here $d$ is either $1$ or $2$. In particular, the boundary
projection multiplicity is $m^2-1$ when $F$ is a line segment,
and $m^2-2$ when $F$ is two-dimensional.
\end{corollary}
\begin{proof}
By Theorem~\ref{thm:normal-face}, the space $\mathcal V$ in
Theorem~\ref{thm:boundary-projection} is spanned by those
coordinate vectors indexed by $\omega_j\in F$. On $\mathcal V$ the
compressed matrix $T_{\mathcal V}$ is diagonal with entries
$\overline{\omega_j}$. Thus
\[
 r_{\mathcal V}=\operatorname{rank}_{\mathbb R}
 \begin{pmatrix}1&\cdots&1\\
 \operatorname{Re}\omega_{j_1}&\cdots&\operatorname{Re}\omega_{j_m}\\
 -\operatorname{Im}\omega_{j_1}&\cdots&-\operatorname{Im}\omega_{j_m}
 \end{pmatrix}=d+1.
\]
The formula follows from Theorem~\ref{thm:boundary-projection}.
Since $|\omega_j|=1$ and $0<q<1$, the minimal face $F$
cannot be a single point, so $d\in\{1,2\}$.
\end{proof}

\begin{example}\label{ex:root-projections}
For $A=\operatorname{diag}(1,\omega,\omega^2)$,
$\omega=e^{2\pi i/3}$, Example~\ref{ex:roots-unity} gives
$W_q(A)=\{\xi:|\xi|\le1\}$ when $0<q\le1/2$.
If $0<q<1/2$, then $F$ is the whole eigenvalue triangle,
so $m_q(A,\xi)=3$ and $p_q(A,\xi)=7$ for every $|\xi|=1$.
At $q=1/2$, the three exceptional boundary points have
$m_q(A,\xi)=2$ and $p_q(A,\xi)=3$; at all other boundary
points these multiplicities are $3$ and $7$, respectively.
\end{example}

\begin{corollary}\label{cor:weighted-projections}
Let $A$ be a positive weighted nilpotent shift, with $r$, $k_L$,
and $q_A$ as in Theorem~\ref{thm:weighted}. For every
$\xi\in\partial W_q(A)$ and $0<q<1$,
\[
 p_q(A,\xi)=
 \begin{cases}
 r^2-2,&0<q<q_A,\\
 k_L^2,&q=q_A>0,\\
 1,&q_A<q<1.
 \end{cases}
\]
\end{corollary}
\begin{proof}
For $0<q<q_A$, Theorem~\ref{thm:weighted} shows that the
admissible second vectors span the top right singular space $E$ of $A$.
Thus $\mathcal V=E$ in Theorem~\ref{thm:boundary-projection}; the first
vectors are their images under the isometry $U$.
The compression $T$ in Theorem~\ref{thm:boundary-exact} is the
adjoint of a direct sum of strictly triangular weighted-shift
blocks, at least one of which is nonzero. Hence
$I_E$, $\operatorname{Re}T$, and $\operatorname{Im}T$ are
linearly independent over $\mathbb R$. Thus $r_{\mathcal V}=3$ and
Theorem~\ref{thm:boundary-projection} gives $p_q(A,\xi)=r^2-2$.

At $q=q_A>0$, the admissible second vectors span the $k_L$-dimensional
top eigenspace of the Hermitian part of $T$. The compression of
$T$ to this subspace equals $q_A I$: the skew-Hermitian part
has zero compression, as in the proof of
Theorem~\ref{thm:weighted}. Thus
\[
 r_{\mathcal V}=\dim_{\mathbb R}
 \operatorname{span}_{\mathbb R}
 \{I_{\mathcal V},\operatorname{Re}T_{\mathcal V},
   \operatorname{Im}T_{\mathcal V}\}=1,
\]
and Theorem~\ref{thm:boundary-projection} gives
$p_q(A,\xi)=k_L^2-1+1=k_L^2$. For $q>q_A$, the generating-vector space is
one-dimensional, so $p_q(A,\xi)=1$. Rotational similarity gives
the same multiplicities at every boundary point.
\end{proof}
\begin{example}
For the shift in Example~\ref{ex:runs},
\[
p_q(A,\xi)=
\begin{cases}
14,&0<q<1/2,\\
4,&q=1/2,\\
1,&1/2<q<1.
\end{cases}
\]
\end{example}

\begin{example}[Vector and projection multiplicities differ]
\label{ex:hermitian-projection}
Let $A=\operatorname{diag}(-1,1)$, with $0<q<1$.
Every point on the ellipse $\partial W_q(A)$ has vector
multiplicity $m_q(A,\xi)=2$, but its generating first vectors have
fixed nonzero coordinate moduli and arbitrary relative phase.
Their rank-one projections span a real three-dimensional space,
so $p_q(A,\xi)=3< m_q(A,\xi)^2=4$.
In contrast, Theorem~\ref{thm:interior-projection} gives the
full $n^2$ projection span at all interior points when $n\ge3$.
\end{example}

\section{The endpoint \texorpdfstring{$q=0$}{q=0}: geometry and multiplicity}\label{sec:qzero}

We now impose $x^*y=0$. In contrast to the case $0<q<1$,
a defective $2\times2$ matrix has vector multiplicity one at the
origin, although the origin is interior to its range. Retain the
previous definitions of $W_q(A)$, $\mathcal M_q(A,\xi)$ and
$m_q(A,\xi)$ at $q=0$, and put
\[
  r_A(x)=\sqrt{\|A^*x\|^2-|x^*Ax|^2},\qquad
  R_0(A)=\max_{\|x\|=1}r_A(x).
\]
The fact that $W_0(A)$ is a disk with radius
$\min_{c\in\C}\|A-cI\|$ is classical; see
Li, Mehta, and Rodman~\cite[Section~2]{LiMehtaRodman}
and Li and Nakazato~\cite[(1.4)]{LiNakazato1998}.
Our concern is the multiplicity of its generating vectors.

\begin{proposition}[Radial geometry and generating-vector sets]\label{prop:q0radial}
Let $A\in M_n(\C)$, $n\ge2$. Then
\[
 W_0(A)=\{\xi\in\C:|\xi|\le R_0(A)\},\qquad
 R_0(A)=\min_{c\in\C}\|A-cI\|.
\]
For every $\xi\in W_0(A)$, the generating-vector set is given by
\[
 \mathcal M_0(A,\xi)=
 \begin{cases}
 \{x:\|x\|=1,\ r_A(x)\ge |\xi|\},&n\ge3,\\
 \{x:\|x\|=1,\ r_A(x)=|\xi|\},&n=2.
 \end{cases}
\]
Thus $\mathcal M_0(A,\xi)$ and $m_0(A,\xi)$ depend only on $|\xi|$.
For fixed $x$, replacing $y$ by $e^{i\theta}y$ rotates $x^*Ay$
but leaves the first vector unchanged.
The minimizing scalar $c$ is unique.
\end{proposition}
\begin{proof}
For a fixed first vector $x$, the image of the unit sphere in $x^\perp$
under $y\mapsto x^*Ay$ is the disk of radius $r_A(x)$ when
$\dim x^\perp\ge2$, and its boundary circle when $\dim x^\perp=1$.
This proves the displayed formulas for $\mathcal M_0(A,\xi)$.
An eigenvector of $A^*$ can serve as $x$ in a generating pair for
zero. The pair $(x,e^{i\theta}y)$ shows that $W_0(A)$ is rotationally
invariant. By convexity and compactness, its containing a point of
maximal modulus $R_0(A)$ implies that it contains the entire closed
disk of that radius.

For any scalar $c$, orthogonality gives
\[
|x^*Ay|=|x^*(A-cI)y|\le\|A-cI\|,
\]
hence $R_0(A)\le\min_c\|A-cI\|$. For the opposite inequality,
set $q=0$ and $\theta=0$ in the proof of
Theorem~\ref{thm:support-dual}. The coercivity estimate, rank-one
reduction, and minimax calculation all remain valid, giving
$h_{W_0(A)}(0)=\min_t\|A+tI\|$. The left side is $R_0(A)$.

For uniqueness, set $f(c)=\|A-cI\|^2$.  For distinct $c_1,c_2$,
\[
 f\!\left(\frac{c_1+c_2}{2}\right)
 \le \frac{f(c_1)+f(c_2)}2-\frac{|c_1-c_2|^2}{4},
\]
by the parallelogram identity applied to $(A-c_1I)v$ and
$(A-c_2I)v$, followed by the supremum over unit $v$. Two distinct
minimizers would have a midpoint with strictly smaller $f$. A minimizer
exists because $f(c)\to\infty$ as $|c|\to\infty$.
\end{proof}

\begin{theorem}[Complete interior classification at $q=0$]\label{thm:q0interior}
Let $A\in M_n(\C)$ be nonscalar and $R=R_0(A)>0$.
If $n\ge3$, then
\[
 m_0(A,\xi)=n\quad (|\xi|<R),\qquad
 \mathcal M_0(A,0)=\{x:\|x\|=1\}.
\]
If $n=2$, then
\[
 m_0(A,\xi)=2\quad (0<|\xi|<R),
 \qquad
 m_0(A,0)=
 \begin{cases}
 2,&A\text{ is diagonalizable},\\
 1,&A\text{ is defective}.
 \end{cases}
\]
A scalar matrix has $W_0(A)=\{0\}$ and $m_0(A,0)=n$.
\end{theorem}
\begin{proof}
For $n\ge3$ and $|\xi|<R$, choose a maximizing unit vector $x_0$
for $r_A$. Since $r_A(x_0)=R> |\xi|$, continuity and the preceding
proposition put an open neighborhood of $x_0$ on the unit sphere
inside $\mathcal M_0(A,\xi)$. Its span is $\C^n$.
For $\xi=0$, every unit $x$ is admissible as a first vector: choose
$y$ of unit length in $\{x,A^*x\}^\perp$, whose dimension is at
least $n-2\ge1$.

For $n=2$, Proposition~\ref{prop:q0radial} gives
$x\in\mathcal M_0(A,0)$ if and only if $r_A(x)=0$, or equivalently
$A^*x\in\C x$. The span is therefore the span of the eigenvectors
of $A^*$. Its dimension is two if $A$ is diagonalizable and one if
$A$ is defective. The case $0<|\xi|<R$ is proved in the next theorem.
For scalar $A$, every unit $x$ generates zero.
\end{proof}

\begin{theorem}[Complete two-dimensional classification]\label{thm:q02by2}
Let $A\in M_2(\C)$ be nonscalar and put
\[
 R=\left\|A-\frac{\operatorname{tr}A}{2}I\right\|.
\]
Then $W_0(A)=\{\xi:|\xi|\le R\}$ and, for $|\xi|\le R$,
\[
 m_0(A,\xi)=
 \begin{cases}
 1,&\xi=0\text{ and }A\text{ is defective},\\
 1,&|\xi|=R\text{ and }A\text{ is nonnormal},\\
 2,&\text{otherwise}.
 \end{cases}
\]
\end{theorem}
\begin{proof}
Write $B=A-\frac12\operatorname{tr}(A)I$. The eigenvalues of $B$
are $\lambda,-\lambda$, so convexity of the classical numerical range
shows $0\in W(B)$. Choose a unit $u_0$ with $u_0^*Bu_0=0$ and
complete it to an orthonormal basis. Since $\operatorname{tr}B=0$,
in this basis $B$ has the form
\[
 B=\begin{pmatrix}0&a\\b&0\end{pmatrix},
 \qquad (a,b)\ne(0,0).
\]
Given a unit vector $x=(u,v)^\top$, the unit vector
$y=(-\bar v,\bar u)^\top$ is orthogonal to $x$, and
\[
 x^*By=a\bar u^{\,2}-b\bar v^{\,2},\qquad
 r_B(x)=|a\bar u^{\,2}-b\bar v^{\,2}|.
\]
The bound $r_B(x)\le |a||u|^2+|b||v|^2\le\max\{|a|,|b|\}$
is attained at a coordinate vector. Thus $R=\|B\|$, and the
optimal center is $\operatorname{tr}(A)/2$.

If $|a|>|b|$, equality in the radius bound requires $v=0$;
if $|b|>|a|$, it requires $u=0$. When $|a|=|b|$, both coordinate
vectors attain $R$. The boundary multiplicity is therefore two
exactly when $|a|=|b|$, which is the condition for $B$ (and $A$)
to be normal.

At the origin use $a\bar u^2=b\bar v^2$. If $ab\ne0$, a solution
with both coordinates nonzero and the solution obtained by changing
the sign of $v$ are independent. If $ab=0$, exactly one of $a,b$
is zero and all solutions lie on one coordinate line. Since
$B^2=abI$, the former case gives two distinct eigenvalues and the
latter gives a nonzero nilpotent $B$.

Finally, let $0<\rho<R$. We construct a unit vector
$x=(u,v)^\top$ with $uv\ne0$ and $r_B(x)=\rho$.
For $ab\ne0$, select the phases so that $a\bar u^2$ and
$b\bar v^2$ have equal arguments. With $t=|u|^2$ we have
\[
 r_B(x)=\bigl||a|t-|b|(1-t)\bigr|.
\]
On $[0,1]$ this continuous function takes the value zero at
$t=|b|/(|a|+|b|)$ and the value $R$ at one endpoint.
Consequently $r_B(x)=\rho$ for some $0<t<1$.
When $ab=0$, take $|u|^2=\rho/|a|$ if $a\ne0$, or
$|v|^2=\rho/|b|$ if $b\ne0$. Both choices give nonzero coordinates.
Changing $v$ to $-v$ leaves $r_B(x)$ unchanged and gives a second,
independent vector. Multiplying the corresponding second vectors by
suitable phases gives any prescribed $\xi$ with $|\xi|=\rho$.
Hence $m_0(A,\xi)=2$ for $0<|\xi|<R$.
\end{proof}

\begin{theorem}[Boundary compression at $q=0$]\label{thm:q0boundary}
Let $A\in M_n(\C)$ be nonscalar, and let $c$ be the unique scalar minimizing $\|A-cI\|$. Set
\[
B=A-cI,\quad R=\|B\|,\quad E=\ker(B^*B-R^2I),
\quad U=R^{-1}B|_E,\quad T=P_EU^*|_E.
\]
Then $0\in W(T)$ and, for $|\xi|=R$,
\[
\mathcal M_0(A,\xi)=\{Uv:v\in E,\ \|v\|=1,\ v^*Tv=0\},
\]
\[
m_0(A,\xi)=\dim_{\C}\spn\{v\in E:\|v\|=1,\ v^*Tv=0\}.
\]
In particular, $1\le m_0(A,\xi)\le\dim E$, and $m_0(A,\xi)=\dim E$ if and only if $0\in\operatorname{ri}W(T)$.
\end{theorem}
\begin{proof}
Write $\xi=Re^{i\theta}$ and let $(x,y)$ generate $\xi$. Because
$x^*y=0$, we have $x^*By=Re^{i\theta}$ and
\[
 R=|x^*By|\le\|By\|\le\|B\|=R.
\]
Hence $y\in E$ and $By=Re^{i\theta}x$. Setting $v=e^{-i\theta}y$
gives $x=Uv$ and
\[
 0=x^*y=e^{i\theta}v^*U^*v=e^{i\theta}v^*Tv.
\]
Conversely these identities show that every unit $v\in E$ with
$v^*Tv=0$ yields the pair $(Uv,e^{i\theta}v)$ generating $\xi$.
This proves the set formula and $0\in W(T)$. Since $U$ is isometric,
taking spans yields the dimension formula; the fixed-expectation
lemma gives the criterion for dimension $\dim E$.
\end{proof}

\begin{corollary}[Normal matrices and the minimal enclosing circle]
\label{cor:q0normal}
Let $A$ be nonscalar and normal, with eigenvalues
$\lambda_1,\dots,\lambda_n$, counted with multiplicity. Let $c$ be
the center and $R$ the radius of the smallest closed disk containing
the eigenvalues, and set
\[
 S=\{j:|\lambda_j-c|=R\},\qquad
 \nu_j=(\lambda_j-c)/R\quad(j\in S),\qquad
 P=\operatorname{conv}\{\nu_j:j\in S\}.
\]
Let $F$ be the unique smallest face of the convex polygon $P$ whose
relative interior contains $0$. Then every boundary point satisfies
\[
 m_0(A,\xi)=\#\{j\in S:\nu_j\in F\},\qquad |\xi|=R.
\]
In particular, $m_0(A,\xi)=|S|$ if and only if $0\in\operatorname{ri}P$.
\end{corollary}
\begin{proof}
For a normal matrix,
$\|A-cI\|=\max_j|\lambda_j-c|$, hence $c$ is the unique center of
the smallest enclosing disk. In a diagonalizing basis, the top
singular space $E$ is the coordinate space indexed by $S$, and the
compression $T$ from Theorem~\ref{thm:q0boundary} is diagonal with
entries $\bar\nu_j$. Thus its zero-numerical vectors are those
$v\in E$ satisfying
$\sum_{j\in S}|v_j|^2\nu_j=0$.
The existence of such a unit vector implies $0\in P$. A positive
coefficient $|v_j|^2$ is possible only if $\nu_j\in F$ by the face
property. Conversely, $0\in\operatorname{ri}F$ permits a convex
representation with positive coefficients on all the indexed points
of $F$, by the averaging argument in the proof of
Theorem~\ref{thm:normal-face}. Independent changes of coordinate signs
then show that the span contains every corresponding coordinate axis.
Taking dimensions proves the formula.
\end{proof}

\begin{example}[Active eigenvalues and boundary multiplicity]
Take $A=\operatorname{diag}(-1,1,i)$. The optimal center is $c=0$,
$R=1$, and the top singular space is all of $\C^3$. However,
$0$ is on the edge joining $-1$ and $1$ in the triangle
$\operatorname{conv}\{-1,1,i\}$. By
Corollary~\ref{cor:q0normal}, the smallest face containing $0$
is this edge. Thus $m_0(A,\xi)=2$ for every $|\xi|=1$,
although $\dim E=3$.
\end{example}

\begin{corollary}[Positive weighted shifts at the endpoint]\label{cor:q0shift}
Let $A$ be a positive weighted nilpotent shift as in
Section~\ref{sec:weighted}, with weights $a_1,\dots,a_{n-1}>0$.
Put $M=\max_j a_j$ and $S=\{j:a_j=M\}$. Then
\[
 W_0(A)=\{\xi:|\xi|\le M\},\qquad
 m_0(A,\xi)=
 \begin{cases}
 n,& |\xi|<M,\quad n\ge3,\\
 |S|,&|\xi|=M.
 \end{cases}
\]
For $n=2$, the multiplicity is one at $\xi=0$ and at $|\xi|=M$,
and two at $0<|\xi|<M$.
\end{corollary}
\begin{proof}
The norm of $A$ equals $M$. Every pair
$(x,y)=(e_j,e_{j+1})$ with $j\in S$ is admissible and generates
$M$. Thus $R_0(A)=M$. At a boundary point,
$|x^*Ay|=M$ forces $y$ into the top right singular space
$\operatorname{span}\{e_{j+1}:j\in S\}$ and
$x$ into $\operatorname{span}\{e_j:j\in S\}$. Conversely, each
coordinate vector $e_j$ with $j\in S$ generates every point of the
boundary circle (rotate $y$). Hence the boundary multiplicity is
$|S|$. The interior assertions follow from
Theorems~\ref{thm:q0interior} and~\ref{thm:q02by2}.
\end{proof}

\begin{example}[A genuine jump at $q=0$]
For the positive weighted shift with weights $(2,1,2)$,
$S=\{1,3\}$, so the boundary multiplicity at zero is $m_0=2$.
The two maximal-weight runs are both singletons; hence $L=1$ and
$q_A=\cos(\pi/2)=0$. Theorem~\ref{thm:weighted} gives $m_q=1$
for every $0<q<1$. Thus the right-hand limit of the boundary
multiplicity at $q=0$ is $1$, not $2$.
\end{example}

\section{Complex parameters in the unit disk}\label{sec:complexq}
Let $q=\rho e^{i\phi}$ with $0<\rho\le1$. The substitution
$y\mapsto e^{-i\phi}y$ changes the inner product from $q$ to $\rho$
and leaves $x$ fixed. The range identity is known from
\cite[Proposition~2.2(c)]{LiMehtaRodman}; here we use the same
substitution at the level of generating vectors and their projections.
There are no admissible pairs if $|q|>1$, by Cauchy--Schwarz.
The case $q=0$ was considered in Section~\ref{sec:qzero}.

\begin{theorem}\label{thm:complex-reduction}
Let $A\in M_n(\C)$, $n\ge2$, and let
$q=\rho e^{i\phi}$ with $0<\rho\le1$. Then
\begin{align}
W_q(A)&=e^{i\phi}W_\rho(A),\label{eq:complex-range}\\
\mathcal M_q(A,\xi)&=\mathcal M_\rho(A,e^{-i\phi}\xi),
\quad \xi\in W_q(A),\label{eq:complex-fiber}\\
m_q(A,\xi)&=m_\rho(A,e^{-i\phi}\xi),\label{eq:complex-mult}\\
p_q(A,\xi)&=p_\rho(A,e^{-i\phi}\xi).\label{eq:complex-projection}
\end{align}
In particular, for $0<|q|<1$,
$m_q(A,\xi)=n$ whenever $\xi\in\operatorname{ri}W_q(A)$.
For nonscalar $A$ and $0<|q|<1$, each supporting line of $W_q(A)$
touches the range at exactly one point.
\end{theorem}
\begin{proof}
Let $(x,y)$ be an admissible pair for the parameter $q$ and set
$y'=e^{-i\phi}y$. Then
\[
\|y'\|=1,\qquad x^*y'=\rho,\qquad
x^*Ay'=e^{-i\phi}x^*Ay.
\]
The reverse correspondence is $y=e^{i\phi}y'$. In particular,
it preserves the first vector $x$, so the generating-vector sets
are identical after the indicated rotation, not merely equinumerous.
This proves both the range and the generating-vector-set identities.
It also preserves the projections $xx^*$, so both multiplicities
are unchanged after the indicated rotation. The consequences
follow from Theorems~\ref{thm:all} and~\ref{thm:support-dual}.
\end{proof}

\begin{theorem}[Support and boundary compression for complex $q$]
\label{thm:complex-support}
Let $A\in M_n(\C)$ be nonscalar and $0<|q|<1$.
For $\theta\in\mathbb R$, the support function satisfies
\begin{equation}\label{eq:complex-support}
h_{W_q(A)}(\theta)
=\min_{t\in\C}
\bigl(\|e^{-i\theta}A+tI\|-\operatorname{Re}(qt)\bigr),
\end{equation}
with an attained and unique minimum. For its minimizer $t$, set
\[
C=e^{-i\theta}A+tI,\quad \sigma=\|C\|,\quad
E=\ker(C^*C-\sigma^2I),\quad
U=\sigma^{-1}C|_E,\quad T=P_EU^*|_E.
\]
Then the supported point and its multiplicity are
\[
\xi_\theta=e^{i\theta}(\sigma-qt),\qquad
m_q(A,\xi_\theta)
=\dim_\C\spn\{v\in E:\ \|v\|=1,\ v^*Tv=q\}.
\]
Moreover, $q\in W(T)$,
$1\le m_q(A,\xi_\theta)\le\dim E$, and
$m_q(A,\xi_\theta)=\dim E$ if and only if
$q\in\operatorname{ri}W(T)$.
\end{theorem}
\begin{proof}
Write $q=\rho e^{i\phi}$, where $0<\rho<1$. Since
$W_q(A)=e^{i\phi}W_\rho(A)$,
$h_{W_q(A)}(\theta)=h_{W_\rho(A)}(\theta-\phi)$.
Apply Theorem~\ref{thm:support-dual} in direction
$e^{i(\theta-\phi)}$, denote its scalar shift by $\tau$,
and put $t=e^{-i\phi}\tau$. Because multiplication of a matrix by a
unimodular scalar leaves its norm unchanged,
this yields~\eqref{eq:complex-support}. Multiplication by
$e^{-i\phi}$ carries the unique real-parameter minimizing shift
$\tau$ to the unique complex-parameter minimizing shift $t$.
In addition, the objective
is coercive:
\[
\|e^{-i\theta}A+tI\|-\operatorname{Re}(qt)
\ge(1-|q|)|t|-\|A\|,
\]
so a minimizer exists.

If $\xi_\theta=x^*Ay$ is a supported value, then
\[
\sigma=
\operatorname{Re}(e^{-i\theta}\xi_\theta+qt)
=\operatorname{Re}(x^*Cy)\le |x^*Cy|\le\sigma.
\]
The equalities imply $Cy=\sigma x$ and $y\in E$. Thus $x=Uy$
and the inner product condition reads $y^*Ty=q$. In the reverse
direction, if $v\in E$ is unit and $v^*Tv=q$, the pair $(Uv,v)$
generates $e^{i\theta}(\sigma-qt)$. Taking spans and using the
isometry $U$ gives the multiplicity formula. Apply the
fixed-expectation lemma to $T$ and $q$ for the last assertion.
\end{proof}

\begin{remark}[The two endpoints]
At $q=0$ the phase $\phi$ is irrelevant, and the complete
classification is given in Section~\ref{sec:qzero}; in particular,
there is a two-dimensional interior exception at the origin.
When $|q|=1$, equality in Cauchy--Schwarz forces $y=qx$, so
\[
W_q(A)=qW(A),\qquad
\mathcal M_q(A,\xi)=\{x:\ \|x\|=1,\ x^*Ax=\xi/q\}.
\]
Thus the classical numerical-range results of Embry recalled in the
Introduction apply directly. Flat boundary segments may occur,
so the strict convexity statement for $|q|<1$ does not extend to
the unit circle. 
\end{remark}

\begin{remark}
One may also consider the complex span of the rank-one matrices
$yx^*$:
\[
d_q(A,\xi)=\dim_{\C}\spn_{\C}
\{yx^*: \|x\|=\|y\|=1,\ x^*y=q,\ x^*Ay=\xi\}.
\]
These matrices belong to $M_n(\C)$ and satisfy
$\tr(yx^*)=q$ and $\tr(Ayx^*)=\xi$.
Thus $d_q(A,\xi)\le n^2-1$ when $q\ne0$ and $A$ is nonscalar.
For $0<q<1$ and $\xi\in\partial W_q(A)$, the boundary
compression in Theorem~\ref{thm:boundary-exact} gives
$yx^*=yy^*U^*$ and $xx^*=Uyy^*U^*$.
As the real and complex spans of a family of Hermitian matrices
have equal respective dimensions, $d_q(A,\xi)=p_q(A,\xi)$
at such boundary points.
We leave a systematic study of $d_q$, particularly at interior
points and at $q=0$, for future work.
\end{remark}

\section*{Acknowledgment}
This research was conducted during Yayi Yuan's six-month visit in academic year 2026--2027 to
the Department of Mathematics and Statistics at the University of
Nevada, Reno.

\section*{Statements and Declarations}

\begin{itemize}
	\item \textbf{Funding}
	
	The research was supported by the Student Science and Technology Innovation Project for the Inner Mongolia Normal University (Grant No. 2025XSKC13) and the Grants for Graduate Students' Research Innovation Projects of Inner Mongolia Normal University (Grant No. KYCXB2025008).
	
	\item \textbf{Conflict of interest/Competing interests}
	
	The authors have no competing interests to declare that are relevant to the content of this article.
	
	\item \textbf{Ethics approval and consent to participate}
	
	Not applicable.
	
	\item \textbf{Consent for publication}
	
	Not applicable.
	
	\item \textbf{Data availability}
	
	Data sharing not applicable to this article as no datasets were generated or analysed during the current study.
	
	\item \textbf{Materials availability}
	
	Not applicable.
	
	\item \textbf{Code availability}
	
	Not applicable.
	
	\item \textbf{Author contribution}
	
	All authors read and approved the final manuscript. All authors contributed equally to this work.
\end{itemize}

\section*{AI disclosure}
The authors confirm that the originality and accuracy of the content have been verified. ChatGPT (version GPT-6) was used to assist with language editing, formatting, and improving the readability and presentation of the manuscript. The authors have checked the terms of use for this AI tool and confirm that its use is suitable for publication. The authors take full responsibility for the integrity of the entire content, including the accuracy of all mathematical statements, proofs, and references.

\bibliographystyle{tfnlm}
\bibliography{references}

\end{document}